\documentclass[9pt]{article}
\usepackage{latexsym,amsfonts,euscript,amsthm}
\usepackage[table]{xcolor}
\usepackage{multirow}
\usepackage{amssymb}
\usepackage{stmaryrd}
\usepackage{mathrsfs,amsmath}
\usepackage{leftidx}
\usepackage[figuresright]{rotating}

\usepackage{booktabs}
\usepackage{tabularx}

\usepackage{tocbibind} 
\usepackage[colorlinks,pagebackref]{hyperref}

\newtheorem{lem}{\bf Lemma}[section]

\newtheorem{prop}[lem]{\bf Proposition}

\newtheorem{thm}[lem]{\bf Theorem}
\newtheorem{rmk}[lem]{\bf Remark}
\newtheorem*{rmkA}{\bf Remark}

\newtheorem{hy}[lem]{\bf Hypothesis}

\newtheorem{mainthm}{Theorem}

\newcommand{\ackname}{Acknowledgements}
\makeatletter
\if@titlepage
  \newenvironment{acknowledgement}{%
    \titlepage
    \null\vfil
    \@beginparpenalty\@lowpenalty
    \begin{center}%
      \bfseries \ackname
      \@endparpenalty\@M
    \end{center}}%
  {\par\vfil\null\endtitlepage}
\else
  \newenvironment{acknowledgement}{%
    \if@twocolumn
      \section*{\ackname}%
    \else
      \small
      \begin{center}%
        {\bfseries \ackname\vspace{-0.5em}\vspace{\z@}}%
      \end{center}%
      \quotation
    \fi}
    {\if@twocolumn\else\endquotation\fi}
\fi
\makeatother

\makeatletter 
\@addtoreset{equation}{section}
\makeatother  

\makeatletter
\def\thanks#1{\protected@xdef\@thanks{\@thanks\protect\footnotetext{#1}}}

	\AtEndDocument{\bigskip{\footnotesize
	\textsc{School of Mathematics and Statistics, Southwest University, No.2 Tian Sheng
	Road,
	Chongqing, 400715, China.} \par  
	\textit{Email address}: \href{mailto:fuming.jiang@hotmail.com}{fuming.jiang@hotmail.com} \par
	\addvspace{\medskipamount}
	\textsc{Department of Mathematics, Suzhou University of Technology, No.99 South Third Ring
	Road,
	Changshu, Jiangsu, 215500, China.} \par  
	\textit{Email address}: \href{mailto:yuzeng2004@163.com}{yuzeng2004@163.com} 
  }}

\title{Finite groups with exactly two nonlinear irreducible $p$-Brauer characters
\thanks{\textbf{Keywords}\,\, character theory, finite groups, $p$-Brauer characters.\\
\textbf{2020 MR Subject Classification}\,\, Primary 20C15, 20C20.\\
The first author is supported by the NSF of China (No. 12471019) and Natural Science Foundation Project of CQ CSTB 
(No. CSTB2024NSCQ-MSX0544).
The second author is supported by the NSF of China (Nos. 12171058, 12326356, 12326358, 12301018), the NSF of Jiangsu Province (No. BK20231356) and the Natural Science Foundation of the Jiangsu Province Higher Education Institutions of China (No. 23KJB110002).}}

\author{Fuming Jiang, Yu Zeng}

\date{}

\begin{document}
\maketitle

\begin{abstract}
	Let $p$ be a prime.
  We classify the finite groups having exactly two irreducible $p$-Brauer characters 
  of degree larger than one.
  The case, where the finite groups have orders not divisible by $p$, was done by P.
  P\'alfy in 1981.
\end{abstract}

\section{Introduction}

Minimal situations constitute a classical theme in group theory. Not only do they
arise naturally, but also they provide valuable hints in searching for general patterns.

For a prime $p$, classifying finite groups with a given relatively small number of irreducible $p$-Brauer characters
is a classical and natural problem in finite group theory.
For instance,
G. Seitz \cite{seitz68} classified the finite $p'$-groups (i.e. finite groups of order not divisible by $p$) having exactly one nonlinear irreducible $p$-Brauer character (i.e. irreducible $p$-Brauer character of degree larger than one);
S. Dolfi and G. Navarro \cite{dolfi} dropped the restriction on the orders of finite groups
and completed the classification started by G. Seitz;
P. P\'alfy \cite{palfy81} gave a complete list of finite $p'$-groups with exactly two nonlinear irreducible $p$-Brauer characters.

The purpose of this paper is to complete the classification of finite groups with exactly two nonlinear irreducible $p$-Brauer characters.
As $\mathrm{O}_{p}(G)$, the \emph{$p$-core} of $G$, is contained in the kernel of every irreducible $p$-Brauer character,
in order to obtain our classification, it is necessary to assume that $\mathrm{O}_{p}(G)=1$.

Before stating the main result of this paper, we introduce some notation.
Let $V$ be an $n$-dimensional vector space over a prime field $\mathbb{F}_\ell$.
As in \cite{manzwolfbook}, we denote by $\Gamma(\ell^n) = \Gamma(V )$ the \emph{semilinear group} of $V$, i.e. (identifying $V$ with $\mathbb{F}_{\ell^n}$)
\[
  \Gamma(\ell^n)=\{ x\mapsto ax^\sigma: x\in \mathbb{F}_{\ell^n}, a\in \mathbb{F}_{\ell^n}^\times, \sigma\in \mathrm{Gal}(\mathbb{F}_{\ell^n}/\mathbb{F}_\ell)\},
\]
and we write $\Gamma_0(\ell^n)=\{ x\mapsto ax: x\in \mathbb{F}_{\ell^n}, a\in \mathbb{F}_{\ell^n}^\times\}$ for the subgroup of multiplications. 
We also denote by $A\Gamma(\ell^{n})=V \rtimes \Gamma(\ell^{n})$ the \emph{affine semilinear group} of $V$.
We use the symbol $\mathsf{PrimitiveGroup}(m,i)$ for the $i$-th primitive group of degree $m$ in the $\mathsf{GAP}$ library of primitive groups (\cite{gap});
$(\mathsf{C}_{m})^{n}$ for the direct product of $n$ copies of cyclic group  $\mathsf{C}_{m}$ of order $m$;
$\mathsf{D}_{2^{n}}$, where $n\geq 3$, for the dihedral group of order $2^{n}$;
$\mathsf{Q}_{2^{n}}$, where $n\geq 3$, for the generalized quaternion group of order $2^{n}$;
$\mathsf{SD}_{2^n}$, where $n\geq 4$, for the semidihedral group of order $2^n$;
$\mathsf{ES}(q^{3}_+)$, where $q$ is an odd prime, for an extraspecial $q$-group of order $q^{3}$ with exponent $q$; 
$\mathsf{Dic}_n$, where $n$ is an odd integer, for the \emph{dicyclic group} of order $4n$, which is isomorphic to $\mathsf{C}_{n}  \rtimes \mathsf{C}_{4}$ where $\mathsf{C}_{4}$ acts by inversion on $\mathsf{C}_{n}$;
$\mathrm{M}_{10}$ for the Mathieu group of degree 10;
$B(2^{n})$ for the (upper triangular) Borel subgroup of $\mathrm{SU}_3(2^{n})$.

\begin{mainthm}\label{thmA}
  Let $G$ be a finite group, and $p$ a prime. 
  Assume that $\mathrm{O}_p(G)=1$.
  Then $G$ has exactly two nonlinear $p$-Brauer characters if and only if either 
  \[
    (G,p)\in \{ (\mathsf{A}_5,5),(\mathsf{S}_5,2), (\mathrm{PGL}_2 (7),2), (\mathrm{M}_{10},2), (\mathrm{Aut}(\mathsf{A}_6 ),2),  (\mathrm{P}\Sigma \mathrm{L}_2 (8),3)\}, 
  \]
  or $G=\mathrm{F}(G) \rtimes H$ is a solvable group where $\mathrm{F}(G)$ is the Fitting subgroup of $G$, and one of the following holds:
  \begin{description}
     \item[(1)] $p\nmid |G|$, and $G$ is either a $2$-group with $|G'|=2$ and $|\mathrm{Z}(G)|=4$,
     or an extraspecial $3$-group;
       \item[(2)] $G$ is a solvable primitive permutation group with a minimal normal subgroup $V=\mathrm{F}(G)\cong (\mathsf{C}_{q})^n$ for a prime $q$ such that $V$ is a primitive $H$-module.
       In this case, either $p=2$ and $G\cong\mathsf{PrimitiveGroup}(q^n,i)$, where
       $$(q^n,i)\in \{(5^2,19),(7^2,25),(3^4,99), (7^2,18),(17^2,82), (23^{2},51)\},$$ 
       or $($up to permutation isomorphism$)$ $H\leq \Gamma(V)=\Gamma(q^n)$
       and one of the following holds:
       \begin{description}
        \item[(2a)] $p\nmid |G|$, and $G\cong \mathrm{PSU}_3(2)$;
        \item[(2b)] $p=2$, and 
        $G$ is isomorphic to one of the following groups: $\mathsf{PrimitiveGroup}(5^2,12)$, $\mathsf{PrimitiveGroup}(3^4,69)$, $A \Gamma(5^2)$ or $A\Gamma(3^4)$; 
        \item[(2c)] $p>2$, $H=L \rtimes C$ where $H'$ is a $p$-group and $C=\mathrm{C}_{H}(v_0)\cong \mathsf{C}_{2}$ for some $v_0\in V$, and $M=V \rtimes L$ is a doubly transitive Frobenius group with kernel $V$ and complement $L$;
        \item[(2d)] $G$ is a Frobenius group with complement $H$ of order $(q^{n}-1)/2$ such that $H'$ is a $p$-group;
        \item[(2e)] $p=n=2$, $H\cong \mathsf{D}_{2(q+1)}\times \mathsf{C}_{(q-1)/2}$ 
        where $q$ is a Mersenne prime;
       \item[(2f)] $p=2$, $H\in \mathrm{Syl}_{2}(\Gamma(q^n))$ where $q^n\in \{ 3^{4}, 5^{2} \}$;
       \item[(2g)] $p=n=2$, $H\cong \mathsf{SD}_{4(q+1)/3}\times \mathsf{C}_{(q-1)/2}$ where $q=3\cdot 2^{k}-1\geq 11$, and $\mathrm{C}_{H}(v_0)=1$ for some $v_0\in V$;
       \end{description} 
       \item[(3)] $G$ is a solvable primitive permutation group with a minimal normal subgroup $V=\mathrm{F}(G)\cong (\mathsf{C}_{q})^n$ for a prime $q$ such that $V=V_1\times V_2$ is an imprimitive $H$-module with $|V_i|=q^k$.
       In this case, $($up to permutation isomorphism$)$ $H\leq T\wr \langle \sigma\rangle$ where $T\leq \Gamma(V_1) =\Gamma(q^{k})$ and $\sigma$ is an involution transposing $V_1$ and $V_2$, and one of the following holds:
       \begin{description}
        \item[(3a)] $p=2$, $q=q^k$ is a Fermat prime, $T=\Gamma_0(q)=\langle t\rangle$ and, either $H= T \wr \langle \sigma\rangle$ or $H=\langle t \sigma\rangle \Phi(T\wr \langle \sigma\rangle)$;
        \item[(3b)] $p=2$, and $G\cong\mathsf{PrimitiveGroup}(3^4,i)$ where $i\in \{36,37, 41, 61,62,63,65,66,67, 82,84,85,86,87,95\}$;
        \item[(3c)] $p>2$, $T=\Gamma_0(q^{k})$ and $H\cong \mathsf{Dic}_{(q^{k}-1)/2}\times \mathsf{C}_{(q^{k}-1)/2}$ where $q^k=2p^s+1$ such that 
        $$(q^{k},p,s)\in \{ (q,p,s), (3^5,11,2), (3^k,(3^k-1)/2,1) \};$$
       \end{description}
       \item[(4)] $p>2$, and $G=VH\times \mathrm{Z}(G)$ where $|\mathrm{Z}(G)|=2$ and 
       $VH$ is a doubly transitive Frobenius group with complement $H$ such that $H'$ is a $p$-group;
          \item[(5)] $(G,p)\in \{ (\mathrm{SL}_2(3),3),(B(8),3), (B(2^{m}),2^{m}+1) \}$
          where $2^{m}+1$ is a Fermat prime larger than $3$;
          \item[(6)] $p=2$, and $\mathrm{F}(G)\cong \mathsf{ES}(q^{3}_+)$ where $q$ is a Mersenne prime.
          Moreover, either $H=\Gamma(q^{2})$ or, $H$ is a maximal subgroup of $\Gamma(q^{2})$, which is isomorphic to either  $\mathsf{C}_{q^{2}-1}$, or $\mathsf{Q}_{2(q+1)}\times \mathsf{C}_{(q-1)/2}$ with $q\geq 7$.
  \end{description}
\end{mainthm}

\begin{rmkA}{\rm
 If (2d) holds, we will see in Proposition \ref{prop: rank3 1} that one of the following holds:
     \begin{description}
      \item[(1)] $H\cong \mathsf{C}_{(q^n-1)/2}$;   
      \item[(2)] $p=n=2$, and $H\cong \mathsf{Q}_{q+1}\times \mathsf{C}_{(q-1)/2}$ where $q$ is a Mersenne prime larger than $3$; 
      \item[(3)] $p$ is a primitive prime divisor of $q^n-1$, $|H|=(q^{n}-1)/2$ and $H'$ is a nontrivial $p$-group.
     \end{description}
     We refer to \cite[Chapter XII, $\S$9]{huppertbook3} for the details of doubly transitive Frobenius groups
that appears in our classification theorem.
  }
\end{rmkA}

Throughout this paper, we only consider finite groups. 
The paper is organized as
follows: in Section 2, we collect auxiliary results; 
in Section 3, we deal with solvable groups with exactly two nonlinear irreducible $p$-Brauer characters;
in Section 4, we classify nonsolvable groups with exactly two nonlinear irreducible $p$-Brauer characters;
in Section 5, we present the proof of Theorem \ref{thmA}.

\section{Auxiliary results}

We mainly follow the notation from \cite{isaacsbook,navarrobook} for character theory and \cite{ccn+85} for
finite simple groups.
Throughout, we consistently refer $p$ and $q$ as primes (except for the subsection 2.1, where $q$ denotes a power of 2).
For a positive integer $n$ and a prime $p$, we write $n_p$ to
denote the maximal $p$-power divisor of $n$.

Let $G$ be a finite group. 
We use $G^\sharp$ to denote the set of nontrivial elements of $G$;
$\pi(G)$ to denote the set of prime divisors of $|G|$;
$\mathrm{Irr}(G)$ to denote the set of irreducible complex characters of $G$;
$\mathrm{Irr}(G)^\sharp$ to denote the set of nontrivial irreducible characters of $G$;
$\mathrm{IBr}_p(G)$ to denote the set of irreducible $p$-Brauer characters of $G$;
$\mathrm{IBr}_{p}(G)^\sharp$ to denote the set of nontrivial irreducible $p$-Brauer characters of $G$;
$\mathrm{IBr}_p^{1}(G)$ to denote the set of nonlinear irreducible $p$-Brauer characters of $G$;
$\omega_{p'}(G)$ to represent the set of orders of $p$-regular elements of $G$.
Let $N\unlhd G$ and $\theta \in \mathrm{IBr}_{p}(N)$.
We identify $\varphi \in \mathrm{IBr}_{p}(G/N)$ (resp. $\chi\in \mathrm{Irr}(G/N)$) with its inflation and view $\mathrm{IBr}_{p}(G/N)$ (resp. $\mathrm{Irr}(G/N)$) as a subset
of $\mathrm{IBr}_{p}(G)$ (resp. $\mathrm{Irr}(G)$).
We also use  $\mathrm{IBr}_p(G|\theta)$ to denote the set of irreducible $p$-Brauer characters of $G$ lying over $\theta$
and $\mathrm{IBr}_{p}(G|N)$ to denote the complement of $\mathrm{IBr}_p(G/N)$ in $\mathrm{IBr}_p(G)$.
Other notation will be recalled or defined when necessary.

We begin with two results in number theory.
 Let $q$ be a prime and $n$ a positive integer.
We call a prime $\ell$ a \emph{primitive prime divisor} for $q^n-1$, if $\ell$ divides $q^n-1$ but does not divide $q^k-1$ for any positive integer $k<n$.
If $\ell$ is a primitive prime divisor for $q^n-1$, then $n$ is the order of $q$ modulo $\ell$.
Hence, $n\mid \ell-1$.
The next lemma is a weak version of Zsigmondy's theorem which is proved by K. Zsigmondy in \cite{zsig1892}, one could also find this result in \cite{huppertbook2}.

\begin{lem}[\mbox{\cite[Chapter IX, Theorem 8.3]{huppertbook2}}]\label{lem: zsig}
	Let $q$ be a prime, and $n$ a positive integer. Then a primitive prime divisor exists for $q^n-1$, except when  $n=2$ and $q$ is a Mersenne prime or $q^n=2^6$. 
\end{lem}

The next lemma follows directly from Lemma \ref{lem: zsig}.

\begin{lem}[\mbox{\cite[Chapter IX, Lemma 2.7]{huppertbook2}}]\label{lem: weakcatalan}
  Let $p$ and $q$ be primes, and $a,b$ be positive integers such that $p^a =
  1 + q^b$.
  Then one of the following cases holds:
\begin{description}
  \item[(1)] $q^b=2^3$ and $p^a = 3^2$;
  \item[(2)] $q=2$, $a=1$, $b$ is a power of $2$ and $p = 2^b+1$ is a Fermat prime;
  \item[(3)] $p=2$, $b=1$, $a$ is a prime and $q = 2^a-1$ is a Mersenne prime.
\end{description}
\end{lem}

\begin{lem}\label{lem: compofpgp}
  Let $G$ be an abelian group, and $P$ be its Sylow $p$-subgroup.
  If $X\cong \mathsf{C}_{p}$ is a maximal cyclic subgroup of $P$,
  then $X$ has a complement in $G$. 
\end{lem}
\begin{proof}
  Let $G$ be a counterexample of minimal possible order.
  Observe that $X$ has a complement in $G$ if and only if it has a complement in $P$,
  and hence it follows by the minimality of $G$ that $G=P$.
  As $X$ is not complemented in $G$ and $X\cong \mathsf{C}_p$, $X\leq \Phi(G)$.
  Write $X=\langle x\rangle$.
  Since $G$ is an abelian $p$-group, $\Phi(G)=\{ g^p:g\in G \}$.
  So, $x=g^p$ for some $g \in G$.
  Consequently, $X< \langle g\rangle$, which contradicts the maximality of $X$.
\end{proof}

\subsection{$p$-groups}

We begin with two elementary results related to $p$-groups.

\begin{lem}\label{lem: wreath2}
  Let $P=A \wr Y$ where $A=\langle a\rangle$ has order $2^n$ and $Y=\langle y\rangle$ has order $2$.
  Then
  \begin{description}
    \item[(1)] $\langle ay\rangle \Phi(P)=\langle ay, a^{2}\rangle$.
    \item[(2)] A maximal subgroup $H$ of $P$ has only three involutions if and only if either $H=A\times A^y$ or $H=\langle ay\rangle \Phi(P)$.
  \end{description}
\end{lem}
\begin{proof}
  As $P$ is a nonabelian $2$-group generated by $a$ and $y$, we deduce that $\Phi(P)$ has index $4$ in $P$.
  Also, due to $A\cap A^y=1$, $A\cap \mathrm{Z}(P)=1$.
  Therefore, we conclude that $\langle a^2\rangle\cap \langle aa^y\rangle=1$. 
  In particular, $\langle a^2, aa^y\rangle=\langle a^2\rangle\times \langle aa^y\rangle$.
   Writing $N=\langle a^2, aa^y\rangle$ and, noting that $P/N\cong (\mathsf{C}_{2})^{2}$, 
   we deduce that $N=\Phi(P)$.
  Finally, it is routine to verify that $\langle ay\rangle N=\langle ay, a^{2}\rangle$.

  Observe that $P/N=\{ N, Ny, Na, Nay \}$, and so $P$ has three maximal subgroups: $N\langle y\rangle$, $N\langle a\rangle$ and $N \langle ay\rangle$.
  Since $N \langle y\rangle$ has at least 4 involutions, either $H=N\langle a\rangle=A\times A^y$ or $H=N \langle ay\rangle$.

  Conversely, suppose either $H=N\langle a\rangle=A\times A^y$ or $H=N \langle ay\rangle$.
  If $H=A\times A^y$, as $A\cong \mathsf{C}_{2^{n}}$, then $H$ has precisely three involutions.
  Now, assume that $H=N \langle ay\rangle$.
  Given that $N\cong \mathsf{C}_{2^{n-1}}\times \mathsf{C}_{2^n}$ with $n\geq 2$, we establish that $N$ contains exactly three involutions.
  To demonstrate that $H$ has only three involutions, it suffices to show that $H-N$ contains no involution.
  For each $h\in H-N$, $h=a^{2i}(aa^y)^kay$ for some integers $i$ and $k$.
  Consequently, a straightforward verification confirms that $h$ cannot be an involution.
\end{proof}

For a $p$-group $P$,
we define the
characteristic subgroup $\Omega_1(P):=\langle x\in P: x^{p}=1\rangle$,
and we denote by $\exp (P)$ the exponent of $P$.

\begin{lem}\label{lem: omega}
  Let $P$ be a nonabelian special $p$-group and $A\leq \mathrm{Aut}(P)$.
  Assume that $A$  
  acts transitively on both $(P/\mathrm{Z}(P))^\sharp$ and $\mathrm{Z}(P)^\sharp$.
  Then 
  every element in $P-\mathrm{Z}(P)$ shares the same order.
  Furthermore, the following hold.
  \begin{description}
    \item[(1)] If $p=2$, then $\Omega_1(P)=\mathrm{Z}(P)$. 
    \item[(2)] If $|P|=p^{3}$ with $p>2$, then $P\cong \mathsf{ES}(p^{3}_+)$. 
  \end{description}
\end{lem}
\begin{proof}
  As $P$ is a nonabelian special $p$-group such that $A$ acts transitively on both $(P/\mathrm{Z}(P))^\sharp$ and $\mathrm{Z}(P)^\sharp$,
  both $P/\mathrm{Z}(P)$ and $\mathrm{Z}(P)$ are elementary abelian $p$-groups.
  Let $x,y\in P-\mathrm{Z}(P)$.
  Since $A$ acts transitively on $(P/\mathrm{Z}(P))^\sharp$,
  there exists some $a\in A$ and $z\in \mathrm{Z}(P)$ such that $y^{a}=xz$.
  Consequently, $o(y)=o(y^{a})=o(xz)=o(x)$ where the last equality holds as $o(z)\leq p\leq o(x)$.

  (1) Assume that $p=2$.
  As $P$ is a nonabelian special $2$-group such that $A$ acts transitively on both $(P/\mathrm{Z}(P))^\sharp$ and $\mathrm{Z}(P)^\sharp$,
  $\exp(P)=4$ and $\exp(\mathrm{Z}(P))=2$.
  Therefore, every element in $P-\mathrm{Z}(P)$ has order $4$. 
  Consequently, $\Omega_1(P)=\mathrm{Z}(P)$.

  (2) Note that $P$ is an extraspecial $p$-group of order $p^{3}$.
   Since there is some element in $P-\mathrm{Z}(P)$ has order $p$,
   all elements in $P^\sharp$ have order $p$.
   Consequently, $P\cong\mathsf{ES}(p^{3}_+)$.
\end{proof}

A finite $2$-group $P$ is a \emph{Suzuki $2$-group} if $P$ is nonabelian, possesses
more than one involution, and has an odd-order automorphism $\sigma$ such that $\langle \sigma\rangle$ transitively permutes all involutions of $P$.
According to G. Higman's classification \cite[Theorem 1]{higman1963}, Suzuki $2$-groups are divided into four classes, namely
$A$, $B$, $C$, $D$.
Furthermore, by \cite[Chapter VIII, Theorem 7.9]{huppertbook2}, we have that 
$$P'=\Phi(P)=\mathrm{Z}(P)=\Omega_1(P)=\{ x\in P:x^2=1 \},$$
and that $|P|=|\mathrm{Z}(P)|^2$ if $P$ is of $A$-type and $|P|=|\mathrm{Z}(P)|^{3}$ otherwise. 
For further details of Suzuki 2-groups, we refer to \cite{higman1963} and \cite[Chapter VIII, $\S$7]{huppertbook2}.

Let $P$ be a Suzuki $2$-group of $B$-type and write $|\mathrm{Z}(P)|=q=2^{m}$.
Following Higman's notation \cite[pages 81 and 82]{higman1963},
we denote by $P=B(m,\theta,\varepsilon)$ which is defined on the set
 of $\mathbb{F}_{q}\times \mathbb{F}_{q}\times \mathbb{F}_{q}$
with multiplication
\[
(a,b,c)\cdot (d,e,f)=(a + d, b + e, c + f + ad^{\theta} + be^{\theta} + \varepsilon ae^{\theta})
\]
where $\theta\in \mathrm{Aut}(\mathbb{F}_{q})$ has odd order and $\varepsilon \in \mathbb{F}_{q}^{\times }$ such that $\varepsilon \neq x^{-1}+x^{\theta}$ for all $x\in \mathbb{F}_{q}^\times $.
By calculation, $\mathrm{Z}(P)=\{ (0,0,c): c\in \mathbb{F}_{q}\}$.
A concrete example of Suzuki 2-groups of $B$-type is the (upper triangular) Sylow 2-subgroup $P(q)$ of $\mathrm{SU}_3(q)$, i.e.
\[
 P(q)=\left\{ \begin{pmatrix}
     1& a & b\\
     0& 1 & a^q\\
     0& 0 & 1\\
 \end{pmatrix} : a,b\in \mathbb{F}_{q^{2}},~b+b^q=a^{q+1}\right\}
\]
(see, for instance, \cite[Kapitel II, Satz 10.12]{huppertbook1}).
By calculation, 
\[\mathrm{Z}(P(q))=\left\{ \begin{pmatrix}
  1& 0 & c\\
  0& 1 & 0\\
  0& 0 & 1\\
\end{pmatrix}:c\in \mathbb{F}_q \right\}.\]
Actually, by applying a result of P. Landrock, we are able to show that $B(m,1,\varepsilon)\cong P(q)$ where $q=2^{m}$.

\begin{lem}\label{lem: theta=1}
  $B(m,1,\varepsilon)\cong P(q)$ where $q=2^{m}$.
\end{lem}
\begin{proof}
  Let $P=B(m,1,\varepsilon)$ and $Z=\mathrm{Z}(P)$.
  Then $Z=\{ (0,0,c):c\in \mathbb{F}_q \}$.
  So, $P/Z=\{ (a,b,0)Z:a,b \in \mathbb{F}_q \}$.
  Set $V=P/Z$, $U=\{ (a,0,0)Z: a  \in \mathbb{F}_q \}$ and $W=\{ (0,b,0)Z: b \in \mathbb{F}_q \}$.
  Then $V=U\times W$.
  Define a map $\sigma:P\to P$ by setting $((a,b,c))^{\sigma}=(\lambda a,\lambda b,\lambda^{2}c)$ where $\lambda$ is a generator of the cyclic group $\mathbb{F}_q^\times$.
  It is routine to check that $\sigma\in \mathrm{Aut}(P)$, and that $\langle \sigma\rangle$ acts transitively on $Z^\sharp$, $U^\sharp$ and $W^\sharp$, respectively.
  In particular, $Z$, $U$ and $W$ are faithful irreducible $\langle \sigma\rangle$-modules over $\mathbb{F}_2$.
  Now, we view $Z$, $U$ and $W$ as $\langle \sigma\rangle$-modules over $\mathbb{F}_2$.
  Observe that the set of eigenvalues of $\sigma$ on $Z$ over $\mathbb{F}_q$ is $\{ (\lambda^2)^{\theta}: \theta\in \mathrm{Aut}(\mathbb{F}_q) \}$,
  and that
  the sets of eigenvalues of $\sigma$ on $U$ and on $W$ over $\mathbb{F}_q$ are $\{ \lambda^{\theta}: \theta\in \mathrm{Aut}(\mathbb{F}_q) \}$.
  Therefore, $\sigma$ shares the same set of eigenvalues on $Z$ and on $U,W$ over $\mathbb{F}_q$.
  Consequently, both $U$ and $W$ are isomorphic to $Z$ as an irreducible $\langle \sigma\rangle$-module over $\mathbb{F}_2$.
  Applying \cite[Lemma 2.5]{landrock1975}, we conclude that either $P\cong P(q)$
  or $P$ is isomorphic to a Sylow $2$-subgroup of $\mathrm{PSL}_3(q)$.
  However, as $\mathrm{Z}(Q)<\Omega_1(Q)$ for $Q\in \mathrm{Syl}_{2}(\mathrm{PSL}_3(q))$,
  $P$ is not isomorphic to $Q$.
  Consequently, $P\cong P(q)$.
\end{proof}

\begin{lem}\label{lem: suzuki2}
  Let $P$ be a Suzuki $2$-group having an 
  odd-order automorphism $\sigma$ such that $\langle \sigma\rangle$ transitively permutes all involutions of $P$.
  Assume that $|P|=2^{3m}$ where $|\mathrm{Z}(P)|=2^m$.  
  Then the following hold.
  \begin{description}
    \item[(1)] If $m$ is a power of $2$, then $P$ is of $B$-type.
    \item[(2)] If $m=3$ and $P/\mathrm{Z}(P)$ is a homogeneous $\langle \sigma\rangle$-module, then $P$ is of $B$-type. 
  \end{description}
\end{lem}
\begin{proof}
  Since $|P|=2^{3m}$ where $|\mathrm{Z}(P)|=2^m$, $P$ is one of the following type: $B$-type, $C$-type or $D$-type.
  If $m$ is a power of $2$, 
  we deduce that $P$ must be of $B$-type (see \cite[Pages 81, 82]{higman1963}).
  If $m=3$, then $P$ is of either  $B$-type or $C$-type (see \cite[Pages 81, 82]{higman1963}).
  Also, as $P/\mathrm{Z}(P)$ is a homogeneous $\langle \sigma\rangle$-module, we deduce from \cite[Page 81, Columns III and IV]{higman1963}
  that $P$ is of $B$-type.
\end{proof}

\begin{prop}\label{prop: suzuki2}
  Let $P=B(m,\theta,\varepsilon)$, $q=2^{m}$ and set $Z=\mathrm{Z}(P)$.
  Assume that $P$ possesses an automorphism $\phi$ of order $q^{2}-1$ such that $\langle \phi\rangle$ acts transitively on both $(P/Z)^\sharp$ and $Z^\sharp$.
  If either $q+1$ is a prime or $q=8$,
  then $P$ is isomorphic to $P(q)$ and $P \rtimes \langle \phi\rangle$ is isomorphic to $B(q)$. 
\end{prop}
\begin{proof}
  Let $\sigma=\phi^{q+1}$ and $\psi=\phi^{q-1}$.
  Then $\langle \psi\rangle$ acts trivially on $Z$, and $\langle \sigma\rangle$ transitively permutes
  all involutions of $P$.
  Without loss of generality, we may assume that $((a,b,c))^{\sigma}=(\lambda a,\lambda b,\lambda\lambda^{\theta}c)$ where $\lambda$ is a generator of the cyclic group $\mathbb{F}_q^\times$ (see the proof of \cite[Lemma 12]{higman1963}).
  Note that 
  $Z=\{ (0,0,c): c\in \mathbb{F}_{q} \}$,
  and hence we set $V=P/Z= \{ (a,b,0)Z:a,b\in \mathbb{F}_q \}$,
  $U=\{ (a,0,0)Z: a\in \mathbb{F}_q \}$ and $W=\{ (0,b,0)Z:b\in \mathbb{F}_q \}$.
  Then $U$ is isomorphic to $W$ as an irreducible $\langle \sigma\rangle$-module over $\mathbb{F}_2$.
  As, by \cite[Proposition 8.2]{doerk92},
  the number of irreducible $\langle \sigma\rangle$-submodules of $V$ over $\mathbb{F}_2$
  is $q+1$,
  it follows that 
  $\langle \psi\rangle$ acts regularly on the set of irreducible $\langle \sigma\rangle$-submodules of $V$ over $\mathbb{F}_2$.
  Consequently,
  we write $W=U^{\delta}$ for some $\delta \in \langle \psi\rangle$.
  
  If $q+1$ is a prime, then $q=2^m$ where $m$ is a power of $2$.
  As $\theta$ is an odd-order element of the $2$-group $\mathrm{Aut}(\mathbb{F}_q)\cong \mathsf{C}_{m}$,
   it follows that $\theta=1$.
  So, Lemma \ref{lem: theta=1} yields that $P\cong P(q)$.

  Assume now that $q=8$.
  Then $P=B(3,\theta,\varepsilon)$, $\mathrm{Aut}(\mathbb{F}_8)\cong \mathsf{C}_{3}$, $o(\psi)=9$ and $o(\delta)\in \{ 3,9 \}$.
  If $\theta=1$, then $P\cong P(q)$ by Lemma \ref{lem: theta=1}.
  Note also that $B(m,\theta,\varepsilon)\cong B(m,\theta^{-1},\varepsilon)$ (see \cite[Page 82]{higman1963}),
  and hence we may assume that $x^{\theta}=x^{2}$ for all $x\in \mathbb{F}_8$.
  Recalling that $W=U^{\delta}$, we deduce that $((a,b,0)Z)^{\delta}=(tb,t^{-1}a+\beta b,0)Z$ for all $a,b \in \mathbb{F}_8$
  where $t\in \mathbb{F}_8^{\times}$ and $\beta=\lambda+\lambda^8$ for some $\lambda \in \mathbb{F}_{8^2}^{\times}$ of order $o(\delta)$.
  As $(0,0,c)^\delta=(0,0,c)$ for all $c\in \mathbb{F}_8$ 
  and 
  $((a,b,c)^{\delta})^{2}= ((a,b,c)^{2})^{\delta}$
  for all $a,b,c\in \mathbb{F}_8$,
  it is routine to check that
  $$tb(tb)^{2}+\varepsilon tb(t^{-1}a)^{2}+\varepsilon tb (\beta b)^{2}+t^{-1}a(t^{-1}a)^{2}+t^{-1}a (\beta b)^{2}+\beta b(t^{-1}a)^{2}+\beta b (\beta b)^{2}=aa^{2}+\varepsilon ab^{2}+bb^2$$
  for all $a,b \in \mathbb{F}_8$.
  In particular, if $a=1$ and $b=0$, then we deduce that $t=1$;
  if $b=1$ and $a=0$, then $\varepsilon \beta^{2}=\beta\beta^{2}$; if $a=b=1$, as $\varepsilon \beta^{2}=\beta\beta^{2}$, then $\beta^{2}=\beta$.
  Therefore, $\beta \in \mathbb{F}_2$.
  As $\beta=\lambda+\lambda^8$ for some $\lambda \in \mathbb{F}_{8^2}^{\times}$ of order $o(\delta)$ where $o(\delta)\in \{ 3,9 \}$,
  we deduce that $\beta=1$.
  Recalling that $\varepsilon \beta^{2}=\beta\beta^{2}$, we conclude that $\varepsilon=\beta=1$.
  However, as the polynomial $x^{3}+x+1$ splits over $\mathbb{F}_8$,
  and so there exists some $y\in \mathbb{F}_8^{\times }$ such that $1=y^{2}+y^{-1}=y^{\theta}+y^{-1}$, which contradicts the definition of $\varepsilon$.

We now identify $P$ with $P(q)$.
 As $P=P(q)$ is a Sylow $2$-subgroup of $\mathrm{SU}_3(q)$ where $q=2^{m}$,
 $B(q)=P  \rtimes T$ is a Borel subgroup of $\mathrm{SU}_3(q)$ where $T\cong \mathsf{C}_{q^{2}-1}$ is a maximal torus of $\mathrm{SU}_3(q)$ (in fact, $T=\langle \operatorname{diag}(\mu,\mu^{q-1},\mu^{-q})\rangle$ where $\mu\in \mathbb{F}_{q^{2}}^\times $ has order $q^2-1$).
  Since $\mathrm{C}_{T}(P)=1$, we deduce $T$ is isomorphic to a subgroup of $\mathrm{Aut}(P)$, say $X$.
  Set $A=\mathrm{Aut}(P)$, $B=\mathrm{C}_{A}(P/Z)$ and $C=\mathrm{C}_{A}(Z)$.
  Then $B$ and $C $ are $A$-invariant subgroups.
  Also, as $Z=P'=\Phi(P)=\{ g^{2}:g\in P \}\cong (\mathsf{C}_{2})^{m}$, 
 it is routine to check that $B\leq C$.
  Note that $B$, as the central automorphism group of $P$, has order $|\mathrm{Hom}(P/Z,Z)|$.
   So, $B$ is a normal $2$-subgroup of $A$.
  Set $\overline{A}=A/B$.
  Note that $\overline{X}$ is the unique cyclic subgroup of $\overline{A}$ with order $q^{2}-1$ (see the proof of \cite[Theorem 1.1]{landrock1975}), and hence $\overline{X}=\overline{\langle \phi\rangle}$.
  As $(|\overline{X}|,|B|)=1$, the Schur-Zassenhaus theorem implies that $X=\langle \phi\rangle^\sigma$ for some $\sigma\in B$.
  So, $B(q)\cong P  \rtimes X= P  \rtimes \langle \phi\rangle^\sigma=(P  \rtimes \langle \phi\rangle)^{\sigma}\cong P \rtimes \langle \phi\rangle$.
\end{proof}

\subsection{Solvable primitive groups of small rank}

Let $G$ be a solvable group.
It is known that $G$ is a \emph{primitive permutation group} on a set $\Omega$ with $m$ elements
if and only if $G=V  \rtimes H$ where $V=\mathrm{F}(G)\cong (\mathsf{C}_{q})^{n}$ is the unique minimal normal subgroup of $G$
such that $m=q^{n}$ (see, for instance, \cite[Kapitel II]{huppertbook1}).
In addition, if $H$ acts transitively on $V^\sharp$, then we say $G$ is a \emph{solvable doubly transitive permutation group};
 if $H$ has exactly two orbits on $V^\sharp$,
 then we say $G$ is a \emph{solvable primitive permutation group of rank $3$}.

 To classify solvable groups with exactly two $p$-Brauer characters,
 we need the following two deep theorems.

\begin{thm}\label{thm: classification of solvable 2-transitive permutation groups}
  Let $G$ be a solvable doubly transitive permutation group.
  Write $G=V \rtimes H$, where $V=\mathrm{F}(G)$ and $H$ is a complement of $V$ in $G$.
  Then $($up to permutation isomorphism$)$ either $H\leq\Gamma(q^n)$ where $q$ is a prime, or $G=\mathsf{PrimitiveGroup}(q^n,i)$ with
 \[
  \begin{aligned}
    (q^n, i) \in & \{(3^2, 6),(3^2, 7),(5^2, 15),(5^2, 18),(5^2, 19),(7^2, 25),(7^2, 29), \\
      & (11^2, 39),(11^2, 42),(23^2, 59),(3^4, 71),(3^4, 90),(3^4, 99)\}.
    \end{aligned}
    \]
\end{thm}
\begin{proof}
  This is \cite[Hauptsatz]{huppert1957} combined with group identification via $\mathsf{GAP}$ \cite{gap}.
\end{proof}

\begin{thm}\label{thm: classification of solvable primitive permutation groups of rank 3}
  Let $G$ be a solvable primitive permutation group of rank $3$.
  Write $G=V \rtimes H$, where $V=\mathrm{F}(G)$ and $H$ is a complement of $V$ in $G$.
  Then one of the following holds:
  \begin{description}
    \item[(1)] Up to permutation isomorphism, $H\leq \Gamma(q^n)$ where $q$ is a prime;
    \item[(2)] $G\cong\mathsf{PrimitiveGroup}(q^n,i)$ with 
      \[
        \begin{aligned}
          (q^n, i) \in & \{(3^4,54),(3^4,56),(3^4,77),(3^4,78),(3^4,80),(3^4,81),(3^4,98),(3^4,102),(3^4,104),(3^4,107),\\
          &(7^2,18),(7^2,27),(13^2,59),(17^2,82),(17^2,83),(17^2,90),(19^2,83),(23^2,51),(29^2,101), \\
          &(31^2,118),(47^2,56),(7^4,775),(3^6,350),(3^6,351),(2^6,25),(2^6,26),(2^6,33),(2^6,37),(2^6,40)\};
          \end{aligned}  
      \]
    \item[(3)] $V=V_1\times V_2$ is an imprimitive $H$-module such that $\mathrm{N}_H(V_i)$ acts transitively on $V_i^\sharp$, and $V_1^\sharp\cup V_2^\sharp$ and $V-(V_1\cup V_2)$ are the nontrivial orbits of the action of $H$ on $V$.   
  \end{description}
\end{thm}
\begin{proof}
  This is \cite[Theorem 1.1]{foulser} combined with group identification via $\mathsf{GAP}$ \cite{gap}.
\end{proof}

\subsection{Some results in character theory}

We end this section with several lemmas related to character theory.

\begin{lem}\label{lem: linearext}
  Let $N$ be a normal subgroup of a finite group $G$.
  Let $\lambda \in \mathrm{IBr}_{p}(N)$ and $T=\mathrm{I}_{G}(\lambda)$.
 \begin{description}
  \item[(1)] If $\lambda$ extends to $T$, then $|\mathrm{IBr}_p(G|\lambda)|=|\mathrm{IBr}_p(T/N)|$.
  \item[(2)] If $G$ is split over $N$ and $\lambda(1)=1$, then $\lambda$ extends to $T$.
  In particular, $|\mathrm{IBr}_p(G|\lambda)|=|\mathrm{IBr}_p(T/N)|$.
 \end{description}
 \end{lem}
\begin{proof}
  By Clifford's correspondence, we may assume that $G=T$.

  If $\lambda$ extends to $G$,
  then part (1) follows directly from \cite[Corollary 8.20]{navarrobook}.

  Now, assume that $G=N  \rtimes H$ and that $\lambda(1)=1$.
  Let $K=\ker(\lambda)$ and observe that $K$ is normal in $G$ (as $\lambda$ is $G$-invariant).
  Set $\overline{G}=G/K$.
  Then $\overline{N}$ is a cyclic $p'$-group.
  Consider the action of $\overline{H}$ on $\mathrm{IBr}_p(\overline{N})=\mathrm{Irr}(\overline{N})$ and the action of $\overline{H}$ on $\overline{N}$.
  As $\overline{H}$ fixes each element in $\mathrm{Irr}(\overline{N})=\langle \lambda\rangle$, by \cite[Corollary 6.33]{isaacsbook}, it also fixes each element of $\overline{N}$  i.e. $\overline{G}=\overline{N}\times \overline{H}$.
  An application of \cite[Theorem 8.21]{navarrobook} to $\overline{G}$ reveals that
  $\mathrm{IBr}_p(\overline{G}|\lambda)=\{ \lambda\times \mu:\mu \in \mathrm{IBr}_p(\overline{H}) \}$.
  In particular, $\lambda$ extends to $G$.
  By part (1), $|\mathrm{IBr}_p(G|\lambda)|=|\mathrm{IBr}_p(T/N)|$.
\end{proof}

\begin{lem}\label{lem: intersection of kernels of Brauer characters}
  Let $G$ be a finite group and $p$ a prime. 
  Then $\mathrm{O}_{p}(G)=\bigcap_{\varphi \in \mathrm{IBr}_{p}(G)}\ker(\varphi)$.
\end{lem}
\begin{proof}
  Let $\varphi \in \mathrm{IBr}_{p}(G)$ and $\mathfrak{X}$ an $\overline{\mathbb{F}_p}[G]$-representation affording $\varphi$ where $\overline{\mathbb{F}_p}$ denotes an algebraic closure of $\mathbb{F}_p$.
  As $\ker(\varphi)=\ker(\mathfrak{X})$, our desired result follows directly by \cite[Chapter III, Corollary 2.13]{feit82}.
\end{proof}

\begin{lem}\label{lem: stabdual}
  Let $p$ be a prime and let $G$ act via automorphisms on an abelian group $V$.
  Then $\mathrm{C}_G(v)$ is a $p$-group for every $v \in V^{\sharp}$ if and only if $\mathrm{I}_G(\alpha)$ is a $p$-group for every $\alpha \in \operatorname{Irr}(V)^{\sharp}$.
\end{lem}
\begin{proof}
  Assume that $\mathrm{C}_{G}(v)$ is a $p$-group for each $v\in V^\sharp$.
  Let $\alpha \in \mathrm{Irr}(V)^\sharp$ and $T=\mathrm{I}_{G}(\alpha)$.
  For each cyclic subgroup $X$ of $T$, since $X$ has a fixed point $\alpha\in\mathrm{Irr}(V)^\sharp$,
  it also has a fixed point $v_0\in V^\sharp$ by Brauer's permutation lemma, i.e. $X\leq \mathrm{C}_{G}(v_0)$,
  therefore, $X$ is a $p$-group.
  In particular, $T$ is a $p$-group.
  Similarly, one can prove the converse statement.
\end{proof}

Let a finite group $A$ act via automorphisms on a finite group $G$.
By Brauer's permutation lemma \cite[Section 6, Lemma 1]{brauer41}, the number of $A$-orbits in the set of conjugacy classes of $p$-regular elements in $G$ corresponds to the number of $A$-orbits in $\mathrm{IBr}_p(G)$.
For instance, if $\mathrm{IBr}_p(G)$ has $k$ $A$-orbits, then the set of conjugacy classes of $p$-regular elements of $G$ also has $k$ $A$-orbits,
so $|\omega_{p'}(G)|\leq k$.
Also, if $G$ is a finite group with $|\omega_{p'}(G)|\leq 2$,
then $|\pi(G)|\leq 2$, indicating that $G$ is solvable by Burnside's $p^{a}q^{b}$-theorem.
These observations will be used in the following without further reference.

\begin{lem}\label{lem: transitive action of G on V-0}
  Let $V$ be a normal subgroup of a finite group $G$.
  If $G$ acts transitively on $\mathrm{Irr}(V)^\sharp$, then the following hold.
  \begin{description}
    \item[(1)] $G$ acts transitively on $V^\sharp$, and $V$ is an abelian minimal normal subgroup of $G$.
    \item[(2)] $\{ \mathrm{I}_{G}(\lambda):\lambda\in \mathrm{Irr}(V)^\sharp \}$ and $\{ \mathrm{C}_{G}(v): v\in V^\sharp\}$ are both transitive $G$-sets.
    \item[(3)] If $\mathrm{I}_{G}(\mu)$ is cyclic for some $\mu\in \mathrm{Irr}(V)^\sharp$, then $\{ \mathrm{I}_{G}(\lambda):\lambda\in \mathrm{Irr}(V)^\sharp \}=\{ \mathrm{C}_{G}(v): v\in V^\sharp\}$.
  \end{description}
\end{lem}
\begin{proof}
 As $G$ acts transitively on $\mathrm{Irr}(V)^\sharp$,
 $G$ also acts transitively on $\{ v^V: v\in V^\sharp \}$ by \cite[Corollary 6.33]{isaacsbook}.
 So, every nontrivial element in $V$ shares the same order which is a prime $\ell$.
 In particular, $V$ is an $\ell$-group.
  Also, every character in $\mathrm{Irr}(V)$ is linear,
  which implies that $V$ is abelian.
  Since $G$ acts transitively on $\mathrm{Irr}(V)^\sharp$, it also acts transitively on $V^\sharp$ by \cite[Corollary 6.33]{isaacsbook}.
  Thus, $V$ is an abelian minimal normal subgroup of $G$.
 
  Note that by part (1), both $\{ \mathrm{I}_{G}(\lambda):\lambda\in \mathrm{Irr}(V)^\sharp \}$ and $\{ \mathrm{C}_{G}(v): v\in V^\sharp\}$ are transitive $G$-sets.
  Assume that $\mathrm{I}_{G}(\mu)$ is cyclic for some $\mu\in \mathrm{Irr}(V)^\sharp$.
  By \cite[Theorem 6.32]{isaacsbook}, the cyclic group $\mathrm{I}_{G}(\mu)$ fixes some $v_0\in V^\sharp$, i.e. $\mathrm{I}_{G}(\mu)\leq \mathrm{C}_{G}(v_0)$.
  Comparing the orders of $\mathrm{I}_{G}(\mu)$ and $\mathrm{C}_{G}(v_0)$,
  we conclude that $\mathrm{I}_{G}(\mu)=\mathrm{C}_{G}(v_0)$.
  Consequently, $\{ \mathrm{I}_{G}(\lambda):\lambda\in \mathrm{Irr}(V)^\sharp \}=\{ \mathrm{C}_{G}(v): v\in V^\sharp\}$.
\end{proof}

\section{Solvable groups}
This section provides a comprehensive description of 
finite solvable groups with exactly two nonlinear irreducible
$p$-Brauer characters.

Let $H$ be a subgroup of $G$.
For $\varphi \in \mathrm{IBr}_{p}(G)$, we denote by $\mathrm{IBr}_{p}(\varphi_H)$ the set of irreducible
constituents of $\varphi_H$.
We denote a disjoint union of two sets $A$ and $B$ by $A \sqcup B$.

\begin{lem}\label{lem: ibrG/N}
 Let $G$ be a finite group with $\mathrm{O}_p(G)=1$ and $N$ a normal subgroup of $G$.
 Then the following hold.
 \begin{description}
  \item[(1)] Assume that $G$ is nonabelian.
  Then $\mathrm{IBr}_p^1(G)=\mathrm{IBr}_p^1(G/N)$ if and only if $N=1$.
  \item[(2)] If $\mathrm{IBr}_p^1(G)=\mathrm{IBr}_p^1(G/N)\sqcup \{ \varphi \}$ and $\mathrm{IBr}_p^1(G/N)\neq \varnothing$, then $\mathrm{IBr}_{p}(N)^\sharp=\mathrm{IBr}_{p}(\varphi_N)$ and $N$ is solvable.
  In addition, if $N$ is nilpotent, then $N$ is an abelian minimal normal subgroup of $G$ such that $G$ acts transitively on $N^\sharp$.
   \item[(3)] Assume that $G$ is solvable and that $N= \mathrm{F}(G)<G$.
   Then $\mathrm{IBr}_p^1(G)=\mathrm{IBr}_p^1(G/N)\sqcup \{ \varphi \}$ if and only if
   $G$ is a solvable doubly transitive permutation group with an abelian minimal normal subgroup $N$ such that $|N|>2$ and $\mathrm{C}_{G}(x)/N$ is a $p$-group for each nontrivial $x \in N$.
 \end{description}
\end{lem}
\begin{proof}
 (1) If $N=1$, then apparently $\mathrm{IBr}_p^1(G)=\mathrm{IBr}_p^1(G/N)$.
 Conversely, working by contradiction, we assume that $N$ is nontrivial.
 Since $N$ is not a $p$-group, it has a nontrivial irreducible $p$-Brauer character $\lambda$.
 Given that $\mathrm{IBr}_{p}^{1}(G)=\mathrm{IBr}_{p}^{1}(G/N)$,
 all characters in $\mathrm{IBr}_{p}(G|\lambda)$ are linear.
 This implies that $\lambda$ extends to a linear character $\hat{\lambda}\in \mathrm{IBr}_{p}(G|\lambda)$.
 Applying \cite[Corollary 8.20]{navarrobook} to $G$, we conclude that $\varphi \hat{\lambda}\in \mathrm{IBr}_p(G)$ lies over $\lambda$ for every $\varphi \in \mathrm{IBr}_p^1(G/N)$. 
 So, $\mathrm{IBr}_p^1(G)=\mathrm{IBr}_p^1(G/N)=\varnothing$, i.e.
 every character in $\mathrm{IBr}_p(G)$ is linear.
 Therefore, $G'\leq \bigcap_{\varphi \in \mathrm{IBr}_p(G)}\ker(\varphi)=\mathrm{O}_p(G)=1$ where the first equality holds by Lemma \ref{lem: intersection of kernels of Brauer characters}.
 Consequently, we deduce that $G$ is abelian, a contradiction.

 (2) 
 As $\mathrm{IBr}_p^1(G)=\mathrm{IBr}_p^1(G/N)\sqcup \{ \varphi \}$, $N$ has a nontrivial irreducible $p$-Brauer character (for instance, an irreducible constituent of $\varphi_N$).
 Let $\lambda$ be a nontrivial character in $\mathrm{IBr}_{p}(N)$, we claim that $\lambda$ must be an irreducible constituent of $\varphi_N$.
 Otherwise, $\lambda$ extends to a linear character $\hat{\lambda} \in \mathrm{IBr}_{p}(G|\lambda)$.
 So, by \cite[Corollary 8.20]{navarrobook}, $\psi \hat{\lambda}\in \mathrm{IBr}_p(G)$ lies over $\lambda$ for every $\psi \in \mathrm{IBr}_p^1(G/N)$, whereas $\psi \hat{\lambda} \notin \mathrm{IBr}_{p}^{1}(G/N)\sqcup \{ \varphi \}$.
 Consequently, we conclude that $\mathrm{IBr}_{p}(N)^\sharp=\mathrm{IBr}_{p}(\varphi_N)$.
 Applying Clifford's theorem \cite[Corollary 8.7]{navarrobook}, we observe that $G$ acts transitively on $\mathrm{IBr}_{p}(N)^\sharp$.
 So, $N$ is solvable by the observations preceding Lemma \ref{lem: transitive action of G on V-0}.
 In addition, if $N$ is nilpotent, as $\mathrm{O}_{p}(G)=1$, then $N$ is a $p'$-group and $\mathrm{Irr}(N)=\mathrm{IBr}_{p}(N)$.
 Thus, by part (1) of Lemma \ref{lem: transitive action of G on V-0}, the rest statements hold.

 (3) As $\mathrm{O}_{p}(G)=1$ and $N=\mathrm{F}(G)$, $N$ is a $p'$-group.
 Hence, $\mathrm{IBr}_{p}(N)=\mathrm{Irr}(N)$.
 Also, $\mathrm{C}_{G}(N)\leq N$ because $G$ is solvable with $\mathrm{F}(G)=N$.
 Note that $G$ is nonnilpotent, and hence $N>\mathrm{Z}(G)$.
 In particular, $|N|>2$.
 If $\mathrm{IBr}_{p}^{1}(G/N)=\varnothing$, then we are done by applying \cite[Theorem A]{dolfi}.
 So, we may assume that $\mathrm{IBr}_{p}^{1}(G/N)\neq \varnothing$.

 Now, we assume first that $\mathrm{IBr}_p^1(G)=\mathrm{IBr}_p^1(G/N)\sqcup \{ \varphi \}$.
 Then, by part (2), $G$ is a solvable doubly transitive permutation group with an abelian minimal normal subgroup $N$ such that $|N|>2$.
 Since $G$ acts transitively on $N^\sharp$, it also acts transitively on $\mathrm{Irr}(N)^\sharp$ by \cite[Corollary 6.33]{isaacsbook}.
  Let $\lambda\in\mathrm{Irr}(N)^\sharp=\mathrm{IBr}_{p}(N)^\sharp$ and set $T=\mathrm{I}_{G}(\lambda)$.
  Then Clifford's theorem \cite[Corollary 8.7]{navarrobook} implies that $\psi(1)\geq |G:T|=|N^\sharp|>1$ for each $\psi \in \mathrm{IBr}_{p}(G|\lambda)$.
  As $\mathrm{IBr}_{p}^{1}(G)=\mathrm{IBr}_p^1(G/N)\sqcup \{ \varphi \}$, we have $\mathrm{IBr}_{p}(G|\lambda)=\{ \varphi \}$.
  In particular, $|\mathrm{IBr}_{p}(T/N)|=|\mathrm{IBr}_{p}(G|\lambda)|=1$ by part (2) of Lemma \ref{lem: linearext}.
   So, $T/N$ is a $p$-group.
  Note that $\mathrm{I}_{H}(\lambda)\cong T/N$ is a $p$-group for every $\lambda \in \mathrm{Irr}(N)^\sharp$, and hence it follows by Lemma \ref{lem: stabdual} that $\mathrm{C}_{G}(x)/N\cong \mathrm{C}_{H}(x)$ is also a $p$-group for every $x\in N^\sharp$.

  Conversely, we assume that $G$ is a solvable doubly transitive permutation group with an abelian minimal normal subgroup $N$ such that $|N|>2$ and $\mathrm{C}_{G}(x)/N$ is a $p$-group for each nontrivial $x \in N$.
  Let $\lambda\in \mathrm{IBr}_p(N)^\sharp$ and set $T=\mathrm{I}_{G}(\lambda)$.
  As $\mathrm{IBr}_p(N)=\mathrm{Irr}(N)$, Lemma \ref{lem: stabdual} implies that $T/N$ is a $p$-group.
  So, \cite[Theorems 8.9 and 8.11]{navarrobook} yields that $|\mathrm{IBr}_{p}(G|\lambda)|=|\mathrm{IBr}_{p}(T|\lambda)|=1$.
  Set $\mathrm{IBr}_{p}(G|\lambda)=\{ \varphi \}$.
  Since $G$ acts transitively on $N^\sharp$, it also acts transitively on $\mathrm{Irr}(N)^\sharp=\mathrm{IBr}_{p}(N)^\sharp$ by \cite[Corollary 6.33]{isaacsbook}.
  So, Clifford's theorem \cite[Corollary 8.7]{navarrobook} implies that $\mathrm{IBr}_{p}(G)-\mathrm{IBr}_{p}(G/N)=\{ \varphi \}$ and $\varphi(1)\geq |G:T|=|N^\sharp|>1$.
  Consequently, $\mathrm{IBr}_p^1(G)=\mathrm{IBr}_p^1(G/N)\sqcup \{ \varphi \}$.
\end{proof}

\subsection{Solvable doubly transitive permutation groups}

In this subsection, we deal with solvable doubly transitive permutation groups with exactly two nonlinear irreducible
$p$-Brauer characters. 

Let $V$ be an $n$-dimensional vector space over a prime field $\mathbb{F}_q$.
Recall that $\Gamma(V)=\Gamma(q^{n})$ is the semilinear group of $V$.
Let $\Gamma_0(V)=\Gamma_0(q^{n})$ be the subgroup of multiplications of $\Gamma(V)$.
It is known that $\Gamma_0(V)\cong \mathsf{C}_{q^{n}-1}$ acts fixed-point-freely on $V$ and, $\Gamma(V)/\Gamma_0(V)\cong \mathrm{Gal}(\mathbb{F}_{q^{n}}/\mathbb{F}_q)\cong \mathsf{C}_{n}$.

\begin{hy}\label{hy: 2-transitive}
  Let $G$ be a solvable doubly transitive permutation group with $\mathrm{O}_p(G)=1$.
  Write $G=V \rtimes H$, where $V=\mathrm{F}(G)$ is a minimal normal subgroup of $G$ and $H$ is a complement of $V$ in $G$.
  Assume that $|\mathrm{IBr}_{p}^{1}(G)|=2$.
\end{hy}

Assuming Hypothesis \ref{hy: 2-transitive},
we see that $G$ is nonabelian,
and that $V$ is a $p'$-group such that $2<|V|<|G|$.
In particular, $\mathrm{IBr}_{p}(V)=\mathrm{Irr}(V)$.

\begin{lem}\label{lem: rank 2 0}
  Assume Hypothesis \ref{hy: 2-transitive}.
  Let $\lambda\in \mathrm{Irr}(V)^\sharp$ and set $I=\mathrm{I}_{H}(\lambda)$. 
  If, up to permutation isomorphism, $H\leq \Gamma(V)$, then 
  $I$ is a cyclic
  group of order $2^{1-t} p^s$ where $t=|\mathrm{IBr}_{p}^{1}(G/V)|\leq 1$ and $s\geq 0$.
  Moreover, if $t=0$, then $p>2$ and $H$ is nonabelian.
\end{lem}
\begin{proof}
  Let $H_0=H\cap \Gamma_0(V)$.
  Then $H_0$ and $H/H_0$ are cyclic.
  Note that $H_0$ acts fixed-point-freely on $V$, 
  and so $H_0\cap I=1$.
  Therefore, $I$ is cyclic because $H/H_0$ is cyclic.
  Observe that $H$ acts transitively on $V^\sharp$ and that, by the observations before this lemma, $V$ is a $p'$-group with $2<|V|<|G|$.
  So, $H$ acts transitively on $\mathrm{IBr}_{p}(V)^\sharp=\mathrm{Irr}(V)^\sharp$ by \cite[Corollary 6.33]{isaacsbook}.  
  In particular, $|H:I|=|V|-1>1$ and $\mathrm{IBr}_{p}(G|\lambda)=\mathrm{IBr}_{p}(G)-\mathrm{IBr}_{p}(G/V)$.
  Therefore, we deduce that every character in $\mathrm{IBr}_{p}(G|\lambda)$ is nonlinear by Clifford's theorem \cite[Corollary 8.7]{navarrobook}, i.e. $\mathrm{IBr}_{p}(G|\lambda)=\mathrm{IBr}_{p}^{1}(G)-\mathrm{IBr}_{p}^{1}(G/V)$.
  Thus, we have that 
  $$|\mathrm{IBr}_{p}(I/\mathrm{O}_{p}(I))|=|\mathrm{IBr}_{p}(I)|=|\mathrm{IBr}_{p}(G|\lambda)|=|\mathrm{IBr}_{p}^{1}(G)|-|\mathrm{IBr}_{p}^{1}(G/V)|$$ 
  where the second equality holds by part (2) of Lemma \ref{lem: linearext}.
  As $|\mathrm{IBr}_{p}^{1}(G)|=2$, 
  it forces by part (1) of Lemma \ref{lem: ibrG/N} that $|\mathrm{IBr}_{p}(I/\mathrm{O}_{p}(I))|=2^{1-t}$ where $t=|\mathrm{IBr}_{p}^{1}(G/V)|\leq 1$.
  Recalling that $I$ is cyclic, we conclude that $I \cong\mathsf{C}_{2^{1-t} p^s}$ where $s$ is a nonnegative integer.

   Furthermore, if $t=0$, then we have that $I\cong \mathsf{C}_{2p^s}$ and $p>2$.
   In this case, we also claim that $H$ is nonabelian.
   In fact, otherwise, since $V$ is a faithful irreducible $H$-module, $H$ is a cyclic group acting fixed-point-freely on $V$;
   therefore, $I=1$ which contradicts $|I|=2p^s$.
\end{proof}

The next two results describe the solvable doubly transitive permutation groups 
with exactly two nonlinear irreducible $p$-Brauer characters.

\begin{prop}\label{prop: rank 2 1}
  Assume Hypothesis \ref{hy: 2-transitive}.
  If $|\mathrm{IBr}_{p}^{1}(G/V)|=0$, then $p>2$, $H=L \rtimes C$ where $H'$ is a $p$-group and $C=\mathrm{C}_{H}(v_0)\cong \mathsf{C}_{2}$ for some $v_0\in V^\sharp$, and $M=V \rtimes L$ is a doubly transitive Frobenius group with kernel $V$ and complement $L$.
\end{prop}
\begin{proof}
  Write $|V|=q^n$ where $q$ is a prime.
   As $G$ is a solvable doubly transitive permutation group,
  applying Theorem \ref{thm: classification of solvable 2-transitive permutation groups} and checking via $\mathsf{GAP}$ \cite{gap},
  we have (up to permutation isomorphism) $H\leq \Gamma(q^n)$.

  Let $\lambda \in \mathrm{Irr}(V)^\sharp$ and $I=\mathrm{I}_{H}(\lambda)$.
  Note that $|\mathrm{IBr}_{p}^{1}(G/V)|=0$,
  and so, by Lemma \ref{lem: rank 2 0}, 
  $H$ is nonabelian and $I\cong \mathsf{C}_{2p^{s}}$ where $p>2$.
  An application of part (3) of Lemma \ref{lem: transitive action of G on V-0} yields that $I=\mathrm{C}_{H}(v_0)$ for some $v_0\in V^\sharp$.
  Given that $|\mathrm{IBr}_{p}^{1}(H)|=|\mathrm{IBr}_{p}^{1}(G/V)|=0$,
  it follows that $H/\mathrm{O}_{p}(H)$ is an abelian $p'$-group, implying that $H'$ is a $p$-group.
  Hence, the Schur-Zassenhaus theorem yields that $H=P \rtimes X$ where $P=\mathrm{O}_{p}(H)\in \mathrm{Syl}_{p}(H)$ and $X$ is an abelian Hall $p'$-subgroup of $H$.
  Let $I_2$ be the Sylow $2$-subgroup of $I$.
  Up to a suitable conjugation, we may assume that $X$ contains $I_2$.

  Consider the coprime action of $X$ on $P$.
  Set $H_0=H\cap \Gamma_0(q^n)$ and $P_0=P \cap H_0$, and observe that $P_0\unlhd H$.
  For each subgroup $B$ of $X$,
  since $H/H_0$ is cyclic and $P\unlhd H$, we deduce that
   \[
    [P,B]=[P,B,B]\leq [P\cap H_0,B]=[P_0,B]\leq P_0.
   \]
   Since $H$ is nonabelian, we conclude that $P_0>1$.
    Moreover, $\mathrm{C}_{X}(P_0)=\mathrm{C}_{X}(P)$ is central in $H$, because $X$ is abelian.
    Write $Z_0=\mathrm{C}_{X}(P_0)$.
    Since $\mathrm{C}_{H}(V)=1$, it follows that $\mathrm{Z}(H)\cap I=1$.
    In fact, as $\mathrm{Z}(H)\cap I\unlhd H$ and $H$ acts transitively on $V^\sharp$,
    it follows that $\mathrm{Z}(H)\cap I=\mathrm{C}_{\mathrm{Z}(H)}(v_0)=\mathrm{C}_{\mathrm{Z}(H)}(V)\leq\mathrm{C}_{H}(V)=1$. 
    For $x\in X-Z_0$, we claim that if $\langle x\rangle\cap Z_0=1$, then $\langle x\rangle$ acts fixed-point-freely on $P_0$.
    Indeed, as $\langle x\rangle$ acts coprimely on the abelian group $P_0$, $P_0=\mathrm{C}_{P_0}(y)\times [P_0,y]$ for each nontrivial $y\in \langle x\rangle$;
    since $P_0$ is a cyclic $p$-group, it follows that $\mathrm{C}_{P_0}(y)=1$ for each nontrivial $y\in \langle x\rangle$.

  Let $X_2$ be a Sylow $2$-subgroup of the abelian $p'$-group $X$.
  Then $I_2\leq X_2$.
  We now claim that $I_2$ has a complement $Y$ in $X$.
  To see this, as $|I_2|=2$, by Lemma \ref{lem: compofpgp}, it is enough to show that $I_2$ is a maximal cyclic subgroup of $X_2$.
  Now, let $A$ be a maximal cyclic subgroup of $X_2$ containing $I_2$.
  As $\Omega_1(A)=\{ x \in A:x^2=1 \}=I_2$ and $Z_0\cap I_2\leq \mathrm{Z}(H)\cap I=1$,
  it follows that $Z_0\cap A=1$, indicating that $P_0  \rtimes A$ is a Frobenius group with kernel $P_0$.
  It is noteworthy that the Frobenius group $P_0  \rtimes A$ acts faithfully on an elementary abelian $q$-group $V$ such that the kernel $P_0$ (of $P_0A$) acts fixed-point-freely on $V$.
  Therefore, \cite[Theorem 15.16]{isaacsbook} implies $A \leq \mathrm{C}_{H}(w)$ for some $w \in V^\sharp$.
  As $H$ acts transitively on $V^\sharp$, $\mathrm{C}_{H}(w)$ and $I=\mathrm{C}_{H}(v_0)$ are $H$-conjugate.
  Comparing the $2$-parts of $|I_2|$ and $|A|$, we conclude that $I_2=A$ is indeed a maximal cyclic subgroup of $X_2$.
  
  Thus, $X=Y\times I_2$ and hence $H=(P \rtimes Y)  \rtimes I_2$.
  Set $L=PY$.
  As $|L:L\cap I|=|LI:I|=|H:I|=|V|-1$ where $I=\mathrm{C}_{H}(v_0)$ for some $v_0\in V^\sharp$,
  it follows that $L$ acts transitively on $V^\sharp$.
  So, $L$ is an irreducible subgroup of $\Gamma(q^n)$.
  Now, let $M=V \rtimes L$.
  Then $M\unlhd G$, because $M$ has index 2 in $G$.
   Noting that $L\cap I$ is a $p$-group, 
   we deduce by part (2) of Lemma \ref{lem: linearext} that $|\mathrm{IBr}_{p}(M|\lambda)|=|\mathrm{IBr}_{p}(L\cap I)|=1$.
  Given that $L$ acts transitively on $\mathrm{Irr}(V)^\sharp=\mathrm{IBr}_{p}(V)^\sharp$ and that $|\mathrm{IBr}_{p}^{1}(M/V)|=0$ (as $|\mathrm{IBr}_{p}^{1}(G/V)|=0$), it follows that $|\mathrm{IBr}_{p}^{1}(M)|=|\mathrm{IBr}_{p}(M|\lambda)|=1$.
  Observe that $M\unlhd G$ implies that $\mathrm{O}_{p}(M)\leq\mathrm{O}_{p}(G)=1$,
  and recall that $p>2$.
  We conclude by \cite[Theorem A]{dolfi} that $M$ is a doubly transitive Frobenius group with complement $L$ such that $L'$ is a $p$-group.
  In particular, $L\cap I=1$ and $I\cong \mathsf{C}_{2}$.
  Consequently, $G=M  \rtimes I$ where $I=\mathrm{C}_{H}(v_0)\cong \mathsf{C}_{2}$.
\end{proof}

\begin{prop}\label{prop: rank 2 2}
  Assume Hypothesis \ref{hy: 2-transitive}.
  If $|\mathrm{IBr}_{p}^{1}(G/V)|=1$, then one of the following holds:
  \begin{description}
    \item[(1)] $p=2$, and $G\cong\mathsf{PrimitiveGroup}(q^{n},i)$ where $(q^{n},i)\in \{ (5^2,19),(7^2,25),(3^4,99) \}$;
    \item[(2)] $p\nmid |G|$, and $G\cong \mathrm{PSU}_3(2)$;
    \item[(3)] $p=2$, and 
    $G$ is isomorphic to one of the following groups: $\mathsf{PrimitiveGroup}(5^2,12)$, $\mathsf{PrimitiveGroup}(3^4,69)$, $A \Gamma(5^2)$ or $A\Gamma(3^4)$.
  \end{description}
\end{prop}
\begin{proof}
  Write $|V|=q^n$ where $q$ is a prime.
  As $G$ is a solvable doubly transitive permutation group,
   applying Theorem \ref{thm: classification of solvable 2-transitive permutation groups} and checking via $\mathsf{GAP}$ \cite{gap},
   we have that either (1) holds, or (up to permutation isomorphism) $H\leq \Gamma(q^n)$.
  
   Assume now that $H\leq \Gamma(q^n)$.
   Let $\lambda \in \mathrm{Irr}(V)^\sharp$ and $I=\mathrm{I}_{H}(\lambda)$.
   As $|\mathrm{IBr}_{p}^{1}(G/V)|=1$, 
   it follows by Lemma \ref{lem: rank 2 0} that $I\cong \mathsf{C}_{p^{s}}$.
   By part (3) of Lemma \ref{lem: transitive action of G on V-0}, we know that $I=\mathrm{C}_{H}(v_0)$ for some $v_0\in V^\sharp$.
   Let $H_0=H\cap\Gamma_0(q^n)$ and $X\in \mathrm{Hall}_{p'}(H)$.
   Note that $H$ acts transitively on $V^\sharp$,
   and so part (2) of Lemma \ref{lem: transitive action of G on V-0} implies that $\mathrm{C}_{H}(v)\cong \mathsf{C}_{p^s}$ for each $v \in V^\sharp$.
   Moreover, $\mathrm{C}_{X}(v)=1$ for each $v \in V^\sharp$,
   i.e. $X$ acts fixed-point-freely on $V$.

 Let $\overline{H}=H/\mathrm{O}_{p}(H)$.
 As $\mathrm{O}_{p}(\overline{H})=1$, $\mathrm{F}(\overline{H})$ is a $p'$-group.
 Note also that $|\mathrm{IBr}_{p}^{1}(\overline{H})|=|\mathrm{IBr}_{p}^{1}(H)|=|\mathrm{IBr}_{p}^{1}(G/V)|=1$. 
 It follows by part (3) of Lemma \ref{lem: ibrG/N} that either $\overline{H}$ is nilpotent or $\overline{H}$
 is a solvable doubly transitive permutation group with an abelian minimal normal subgroup $\mathrm{F}(\overline{H})$ such that $|\mathrm{F}(\overline{H})|>2$ and $\mathrm{C}_{\overline{H}}(x)/\mathrm{F}(\overline{H})$ is a $p$-group for each nontrivial $x \in \mathrm{F}(\overline{H})$.

 Let us first assume that $\overline{H}$ is nilpotent.
 In this case, $\overline{H}$ is a $p'$-group.
 Since $|\mathrm{IBr}_{p}^{1}(\overline{H})|=1$, \cite[Theorem A]{dolfi} yields that $\overline{H}$ is an extraspecial $2$-group and $p>2$.
 So, the Schur-Zassenhaus theorem implies that $H=\mathrm{O}_{p}(H)  \rtimes X$ where $X$ is an extraspecial $2$-group. 
 As $X$ acts fixed-point-freely on $V$, we conclude that $X\cong \mathsf{Q}_8$ and $X\cap I=1$.
 Hence, $|H:I|=8p^a$ for some nonnegative integer $a$.
 Since $H$ acts transitively on $V^\sharp$,
 we have $|H:I|=q^n-1$.
 Therefore, we have the equation 
 \begin{equation}\label{nileq}
   \begin{aligned}
     q^n-1=8p^a,
   \end{aligned}
 \end{equation}
 where $q$ is an odd prime.
 Because $H_0$ is cyclic, $H_0\cap X$ is also cyclic, implying that $2\mid|X/(H_0\cap X)|=|H_0X/H_0|\mid |H/H_0|\mid n$.
 By (\ref{nileq}), $(q^{n/2}-1)(q^{n/2}+1)=8p^a$.
 As $(q^{n/2}-1,q^{n/2}+1)=2$, we have two possible cases: $q^{n/2}-1=2$ or $q^{n/2}-1=4$.
 If $q^{n/2}-1=4$, then $q=q^{n/2}=5$ and $p=p^a=3$.
 So, $H\leq \Gamma(5^2)$.
 However, $\Gamma(5^2)$ does not have a subgroup isomorphic to $\mathsf{Q}_8$.
 So, $q^{n/2}-1=2$ and $q^{n/2}+1=4$.
 By calculation, we have that $q=q^{n/2}=3$ and $a=0$.
 Note that $|H:I|_p=p^a=1$ and that $I$ is a $p$-group, and hence $I=\mathrm{C}_{H}(v_0)=\mathrm{O}_{p}(H)\unlhd H$.
 Because $H$ acts transitively on $V^\sharp$, we conclude that $[I,V]=1$.
 Consequently, $\mathrm{O}_{p}(H)=I=1$ provided that $\mathrm{C}_{H}(V)=1$.
 Therefore, $G$ is a doubly transitive Frobenius group isomorphic to $(\mathsf{C}_{3})^{2}  \rtimes \mathsf{Q}_8$, which is also known as $\mathrm{PSU}_3(2)$.
 Thus, (2) holds.

 We next assume that $\overline{H}$
 is a solvable doubly transitive permutation group with an abelian minimal normal subgroup $\mathrm{F}(\overline{H})$ such that $|\mathrm{F}(\overline{H})|>2$ and $\mathrm{C}_{\overline{H}}(x)/\mathrm{F}(\overline{H})$ is a $p$-group for each nontrivial $x \in \mathrm{F}(\overline{H})$.
 Then $\mathrm{F}(\overline{H})$ is an elementary abelian $r$-group for some prime $r\neq p$ and $\mathrm{F}(\overline{H})\leq \overline{X}$.  
 Since $\overline{X}=X\mathrm{O}_{p}(H)/\mathrm{O}_{p}(H)$ is isomorphic to a Frobenius complement $X$,
 every Sylow subgroup of $\overline{X}$ is either cyclic or quaternion.
 As $\mathrm{F}(\overline{H})\leq \overline{X}$, we conclude that $\mathrm{F}(\overline{H})\cong \mathsf{C}_r$ where $r>2$. 
 Furthermore, $\overline{H}\cong A\Gamma(r)$.
 Note that $\overline{X}$ is a subgroup of the Frobenius group $\overline{H}$ containing its kernel $\mathrm{F}(\overline{H})$, and that $\overline{X}$ is isomorphic to a Frobenius complement $X$. 
 We deduce that $\overline{X}=\mathrm{F}(\overline{H})\cong \mathsf{C}_{r}$, because a subgroup of a Frobenius group containing the kernel is either a Frobenius group or the kernel itself. 
 Therefore, $X\cong \overline{X}\cong \mathsf{C}_r$.
 Since $\overline{H}$ is a doubly transitive Frobenius group with kernel $\overline{X}\cong \mathsf{C}_{r}$ where $r>2$,
 we have $r-1=|\overline{H}/\overline{X}|=p^b$ for some positive integer $b$. 
 By Lemma \ref{lem: weakcatalan}, $p=2$ and $r=2^b+1$ is a Fermat prime where $b$ is a power of 2.
 Because $H$ is a $\{ 2,r \}$-group with $|H|_r=r$, we set $|H:I|=2^cr$. 
 Hence, $q^n-1=2^cr$ where $c\geq 1$ and $q$ is an odd prime (as $q\neq p=2$).
 Note that 
 $\overline{H_0}\leq \mathrm{F}(\overline{H})=\overline{X}$,
 and so $2\mid (r-1)\mid |\overline{H}:\overline{X}|\mid |H:H_0|\mid n$.
 Thus,
 we have $(q^{n/2}-1)(q^{n/2}+1)=2^cr$.
 Given that $(q^{n/2}-1,q^{n/2}+1)=2$, we conclude that 
 either $q^{n/2}-1=2^{c-1}$ or $q^{n/2}-1=2r$.

 If $q^{n/2}-1=2r$,
 then $q^{n/2}+1=2^{c-1}$.
 According to Lemma \ref{lem: weakcatalan}, we deduce that $q=q^{n/2}=2^{c-1}-1$ and $r=2^{c-2}-1$ are Mersenne primes.
 Since both $c-1$ and $c-2$ are primes, we have $c=4$, which leads to $q=q^{n/2}=7$ and $r=3$.
 In this case, $V\cong (\mathsf{C}_{7})^{2}$ and $H\leq \Gamma(7^2)$ with $48\mid |H|$.
 Recall that $|\mathrm{IBr}_{p}^{1}(G)|=2$ and $|\mathrm{IBr}_{p}^{1}(H)|=1$,
 and so we conclude a contradiction by checking via $\mathsf{GAP}$  \cite{gap}.

 So, $q^{n/2}-1=2^{c-1}$ and $q^{n/2}+1=2r$.
 Applying Lemma \ref{lem: weakcatalan}, we have that either $q=q^{n/2}=2^{c-1}+1$ and $r=2^{c-2}+1$ are Fermat primes, or $q^{n/2}=3^2$, $c=4$ and $r=5$.
 If the former holds, then $c=3$, $q=q^{n/2}=5$ and $r=3$.
 In this case, we have $H\leq \Gamma(5^2)$ with $24\mid |H|$.
 As $|\mathrm{IBr}_{p}^{1}(G)|=p=2$, we check via $\mathsf{GAP}$ \cite{gap} that $G$ is either $\mathsf{PrimitiveGroup}(5^2,12)$ or $A \Gamma(5^2)$.
 If the latter holds, then $V\cong (\mathsf{C}_{3})^{4}$ and $H\leq \Gamma(3^4)$ with $80\mid |H|$.
 Since $|\mathrm{IBr}_{p}^{1}(G)|=p=2$, one checks via $\mathsf{GAP}$ \cite{gap} that $G$ is either $\mathsf{PrimitiveGroup}(3^4,69)$ or $A \Gamma(3^4)$.
\end{proof}

\subsection{Solvable primitive permutation groups of rank 3}

In this subsection, we deal with solvable primitive permutation groups of rank 3 with exactly two nonlinear irreducible
$p$-Brauer characters.

In the following, we use 
$\mathfrak{M}$ to denote the set of Mersenne primes, 
and $\mathfrak{F}$ to denote the set of Fermat primes.

\begin{lem}\label{lem: rank3hgamma} 
  Let $V$ be an $n$-dimensional vector space over a prime field $\mathbb{F}_q$,
   and $G$ a nonabelian irreducible subgroup of the semilinear group $\Gamma(q^{n})$.
   Assume that $G'$ and $\mathrm{C}_{G}(v)$, for all $v\in V^\sharp$, are $p$-groups where $p\neq q$.
   If $V^\sharp$ is a union of two distinct $G$-orbits,
   then 
   one of the following holds:
   \begin{description}
   \item[(1)] $p=n=2$, $G\cong \mathsf{D}_{2(q+1)}\times \mathsf{C}_{(q-1)/2}$ 
   where $q\in \mathfrak{M}$;
   \item[(2)]  $p=n=2$, $G\cong \mathsf{Q}_{q+1}\times \mathsf{C}_{(q-1)/2}$ acts fixed-point-freely on $V$ where $q\in \mathfrak{M}$ is larger than $3$; 
   \item[(3)] $p$ is a primitive prime divisor of $q^n-1$, $G$ acts fixed-point-freely on $V$ such that $|G|=(q^{n}-1)/2$ and $G'$ is a nontrivial $p$-group;
    \item[(4)] $p=2$, $G\in \mathrm{Syl}_{2}(\Gamma(q^{n}))$ where $q^n\in \{ 3^4,5^2 \}$;
    \item[(5)] $p=n=2$, $G\cong \mathsf{SD}_{4(q+1)/3}\times \mathsf{C}_{(q-1)/2}$ where $q=3\cdot 2^{k}-1\geq 11$, and $\mathrm{C}_{G}(v_0)=1$ for some $v_0\in V^\sharp$.
   \end{description}
\end{lem}
\begin{proof}
  Let $P$ be a Sylow $p$-subgroup of $G$, and observe that $G'$ is a nontrivial $p$-group.
  Thus, $1<G'\leq P$ and therefore $P\unlhd G$.
   So, the Schur-Zassenhaus theorem yields that $G=P  \rtimes H$ where $H\in \mathrm{Hall}_{p'}(G)$ is abelian.
  Since $\mathrm{C}_{G}(v)$ is a $p$-group for each $v\in V^\sharp$, $H$ acts fixed-point-freely on $V$.
   In particular, $H$ is cyclic, and for each $v\in V^\sharp$, $|H|$ divides $|G:\mathrm{C}_{G}(v)|$.
   Now, we take two elements $u,w\in V^\sharp$ belonging to distinct $G$-orbits, and write $A=\mathrm{C}_{G}(u)$, $B=\mathrm{C}_{G}(w)$ and $h=|H|$.
   Then $V^\sharp=u^{G}\sqcup w^{G}$ where $|u^{G}|=|G:A|=p^a h$ and $|w^{G}|=|G:B|=p^b h$.
   Consequently, $q^n-1=p^{a}h+p^{b}h$.
   Without loss of generality, we may assume that $a\geq b$.

   Set $\Gamma=\Gamma(q^{n})$, $\Gamma_0=\Gamma_0(q^{n})$, $G_0=G\cap \Gamma_0$ and $P_0=P\cap G_0$.
   Next, we divide the remaining proof into three claims.
   Note that, if $q^{n}-1$ has a primitive prime divisor $\ell$,
   since $G$ is a nonabelian subgroup of $\Gamma(q^{n})$, we have $\ell>2$.

\textbf{Claim 1.} If $q^n-1$ does not have a primitive prime divisor, then either (1) or (2) holds.

  Since $q^n-1$ does not have a primitive prime divisor,
  Lemma \ref{lem: zsig} implies that either $q^n=2^6$ or $n=2$ and $q \in \mathfrak{M}$.
    If the former holds, then $p^{a}h+p^{b}h=63$.
  We claim first that $h>1$.
  Otherwise, we would have $p^b(p^{a-b}+1)=63$,
  but this equation has no solution.
  Next, we show that $b=0$.
  If not, given that $(p,h)=1$, we would have $p^bh=21$.
  This would imply $p^{a-b}=2$, whereas $p\neq q= 2$.
  Recall that $B=\mathrm{C}_{G}(w)$, where $w\in V^\sharp$, and that $|G:B|=p^bh=h$.
  Therefore, $\mathrm{C}_{G}(w)=B=P$. 
  As $G_0$ acts fixed-point-freely on $V$, $P_0\cap P\leq \mathrm{C}_{G_0}(w)=1$, and so $P\cap G_0=P_0=1$.
  Noting that $G/G_0$ is cyclic, we deduce that $P$ is cyclic.
  Since $H$ acts coprimely on $P$, $[P,H]=[P,H,H]\leq [P\cap G_0,H ]=1$. 
  So, $G=P\times H$ is abelian, a contradiction.

  So, $n=2$ and $q \in \mathfrak{M}$.
  As $q$ is a Mersenne prime, $q\equiv -1~(\mathrm{mod}~4)$,
  and so Sylow $2$-subgroups of $\mathrm{GL}_2(q)$ are isomorphic to $\mathsf{SD}_{4(q+1)}$
  (see, for instance, \cite[Page 142]{carterfong1964}).
  Noting that $\Gamma\leq \mathrm{GL}_2(q)$ and that $|\Gamma|_2=|\mathrm{GL}_2(q)|_2$,
  we deduce by \cite[Lemma 6.5]{manzwolfbook} that $\Gamma=R\times S$ where $R\cong \mathsf{SD}_{4(q+1)}$  and $S\cong \mathsf{C}_{(q-1)/2}$ is of odd order.
  Since $G$ is a nonabelian subgroup of $\Gamma$, $G=P\times H$ where $P$ is a nonabelian Sylow $p$-subgroup of $G$.
  This implies that $p=2$ and $P\leq R$.
  As $|P:P_0|\mid |G:G_0|\mid 2$, it follows that $|P:P_0|=2$. 
  Further, taking into account that $q^2-1=2^bh(2^{a-b}+1)$, we claim that $a=b$.
  Indeed, otherwise $a>b$ which implies that $(q^{2}-1)_2=2^{b}$ and $|B|>|A|$;
  since $|\mathrm{GL}_2(q)|_2=2(q^{2}-1)_2=2^{b+1}=4(q+1)$ and $2^{b+1}\mid 2^{a}\mid |P|\mid |\mathrm{GL}_2(q)|_2$, we deduce that $|P|=4(q+1)$;
  as $P\leq R\cong \mathsf{SD}_{4(q+1)}$, $P=R$;
  recall that $B=\mathrm{C}_{H}(w)$ does not fix $u$,
  and so it follows that $\mathrm{C}_{V}(B)\cong \mathsf{C}_{q}$;
  let $P_1$ be the cyclic subgroup of order $q+1$ of $P_0$, and observe that $P_1$ acts regularly on the set of maximal subgroups of $V\cong (\mathsf{C}_{q})^{2}$;
  so, $u \in \mathrm{C}_{V}(B)^{x}$ for some $x\in P_1$, i.e. $B^{x}\leq \mathrm{C}_{G}(u)=A$, whereas $|B|>|A|$, a contradiction.
   Thus, $a=b$, $|A|=|B|$ and $2|G:B|=q^{2}-1$. 
  As a result, $|P|=(q+1)|B|$ and $|H|=(q-1)/2$.
  Since $V^\sharp$ is a union of two $G$-orbits, we have $G<\Gamma$, which leads to $P<R$.
  If $|A|=|B|>1$, then $P$ is a maximal subgroup of $R$. 
  Notably, $\mathsf{SD}_{4(q+1)}$ has only three maximal subgroups: $\mathsf{C}_{2(q+1)}$, $\mathsf{D}_{2(q+1)}$ and $\mathsf{Q}_{2(q+1)}$. 
  This implies $P\cong \mathsf{D}_{2(q+1)}$ and $G\cong \mathsf{D}_{2(q+1)}\times \mathsf{C}_{(q-1)/2}$. 
  In fact, otherwise, as $P$ is nonabelian, $P\cong \mathsf{Q}_{2(q+1)}$;
  since $\mathrm{Z}(P)$ (the unique subgroup of order $2$ of $P$) is contained in both $A$ and $B$,
  it follows that $[\mathrm{Z}(P),V]=1$ which contradicts $\mathrm{C}_{G}(V)=1$. 
  If $|A|=|B|=1$,
  then $G$ acts fixed-point-freely on $V$ and $|G|=(q^{2}-1)/2$.
  Recall that $G$ is nonabelian, and so $G\cong\mathsf{Q}_{q+1}\times \mathsf{C}_{(q-1)/2}$ where $q \in \mathfrak{M}$ is larger than $3$.

\textbf{Claim 2.} If $q^n-1$ has a primitive prime divisor $\ell$ such that $\ell\mid |G|$, then $\ell=p$ and (3) holds.

Let $L\in \mathrm{Syl}_{\ell}(G)$.
As $\ell$ is a primitive divisor of $q^n-1$ and also $\ell\mid |G|$,
\cite[Corollary 6.6, Lemma 6.4]{manzwolfbook} yields that $P\leq \mathrm{F}(G)=\mathrm{C}_{G}(L)=G_0$.
Let $v\in V^\sharp$, and recall that $\mathrm{C}_{G}(v)\leq P$, 
and thus $\mathrm{C}_{G}(v)\leq \mathrm{C}_{G_0}(v)=1$ i.e. $G$ acts fixed-point-freely on $V$.
In particular, $|G|=(q^{n}-1)/2$.
Since $L\leq \mathrm{C}_{G}(L)=G_0$ by \cite[Lemma 6.4]{manzwolfbook}, we conclude that $p=\ell$.
In fact, if $p \neq \ell$, due to $G'\leq P$, we have $[G,L]\leq G'\cap L=1$, whereas $\mathrm{C}_{G}(L)=G_0<G$.
Therefore, (3) holds because $1<G'\leq P$.

\textbf{Claim 3.} If $q^n-1$ has primitive prime divisors and $\pi(G)$ contains no primitive prime divisor of $q^n-1$, then either (4) or (5) holds.

   Let $\ell$ be a primitive prime divisor of $q^{n}-1$.
   Then $\ell \nmid |G|$.
   Let $L$ be a subgroup of $\Gamma$ of order $\ell$.
   According to \cite[Lemma 6.5]{manzwolfbook}, we have $L\leq \mathrm{C}_{\Gamma}(L) =\Gamma_0$.
   Note that every subgroup of the cyclic group $\Gamma_0$ is characteristic, and so $L$ is normal in $\Gamma$.
   Moreover, as $\ell\nmid |G|$, the Schur-Zassenhaus theorem yields that $LG=L \rtimes G$, and every Hall $\ell'$-subgroup of $LG$ is conjugate to $G$ by some element in $L$.
   Since $G$ is an irreducible subgroup of $\Gamma$, so is $LG$.
   Next, we claim that $LG$ acts transitively on $V^\sharp$.
   In fact, as $\ell\nmid |u^G|$ and $\ell \mid |u^{LG}|$, we have $u^{G} \subset u^{LG}$;
   note that $V^\sharp=u^{G}\sqcup w^{G}$ and that $u^{LG}$ is $G$-invariant, and so $V^\sharp=u^{LG}$.
   Observing that $LG$ acts transitively on $V^\sharp$ and that $a\ge b$ and $\ell>2$,
   we deduce that  $q^{n}-1=p^ah+p^{b}h=\ell p^bh$.
   By calculation, we have that $\ell-1=p^{a-b}$.
   Since $\ell>2$, Lemma \ref{lem: weakcatalan} implies that $p=2$ and $\ell=2^{a-b}+1 \in \mathfrak{F}$.

   Note that $G/\mathrm{C}_{G}(L)$ is isomorphic to a subgroup of $\mathrm{Aut}(L)\cong \mathsf{C}_{\ell-1}$ where $\ell-1=2^{a-b}$.
   As $\mathrm{C}_{G}(L)=G\cap \Gamma_0=G_0=P_0\times H_0$,
   it follows that $G=P \mathrm{C}_{G}(L)=P\times H_0$ and $H=H_0$.
   Also, since $G$ is nonabelian and $H$ is cyclic, $P$ is nonabelian.
   Consequently, $P_0>1$. 
   Indeed, otherwise $P$ is cyclic, a contradiction. 
   As $|P/P_0|\mid |G/G_0|=|G/\mathrm{C}_{G}(L)|\mid 2^{a-b}$ and $|P_{0}|\mid |G:B|_2=2^b$,
   we have $2^a=|G:A|_2\leq |P|=|P_0||P/P_0|\leq 2^a$.
   So, $|P|=2^a$, $|G/G_0|=2^{a-b}$, $|P_0|=2^{b}$ and $A=1$.
    As $\ell$ is a primitive prime divisor of $q^n-1$, we have $n\mid \ell-1=2^{a-b}$.
   As $2^{a-b}=|G/G_0|\mid n$, it follows that $n=2^{a-b}$.
   Since $q\neq p=2$, $q^{2^{a-b}}-1=\ell 2^bh$ where $b>0$.
   Recalling that $\pi(G)$ contains no primitive prime divisor of $q^{2^{a-b}}-1$ and that $h\mid |G|$,
   we conclude that every odd prime divisor of $q^{2^{a-b-1}}+1$ is a primitive prime divisor of $q^{2^{a-b}}-1$.
   Indeed, for an odd prime divisor $r$ of $q^{2^{a-b-1}}+1$, 
   we denote $k$ the order of $q$ modulo $r$;
   by Fermat's little theorem, $k\mid 2^{a-b}$;
   since $r\nmid q^{2^{a-b-1}}-1$, $k=2^{a-b}$ i.e. $r$ is a primitive prime divisor of $q^{2^{a-b}}-1$.  
   Hence,
   $h\mid q^{2^{a-b-1}}-1$ (as $h$ is an odd divisor of $q^{2^{a-b}}-1$) and $\ell\mid q^{2^{a-b-1}}+1$.
   Note that $(q^{2^{a-b-1}}+1,q^{2^{a-b-1}}-1)=2$,
   and so either $q^{2^{a-b-1}}+1=2\ell$ or $q^{2^{a-b-1}}+1=2^{b-1}\ell$.

   If $q^{2^{a-b-1}}+1=2\ell$,
   then $q^{2^{a-b-1}}-1=2^{b-1}h$.
   As $\ell=2^{a-b}+1$, $q^{2^{a-b-1}}=2^{a-b+1}+1$.
   Recall that $n=2^{a-b}$.
   By Lemma \ref{lem: weakcatalan}, either $a-b=2$, $q^{n}=3^{4}$ and $\ell=5$ or $a-b=1$, $q^{n}=5^{2}$ and $\ell=3$.
   If the former holds, then $h=1$, $b=4$ and $a=6$.
   If the latter holds, then $h=1$, $b=3$ and $a=4$.
   In both cases, one has that $G\in \mathrm{Syl}_{2}(\Gamma)$.
   So, (4) holds.

   If $q^{2^{a-b-1}}+1=2^{b-1}\ell$,
   then $q^{2^{a-b-1}}-1=2h$.
   As $h$ is odd, $(2h)_2=2$, and so $a-b=1$.
   By calculation, we have that $\ell=2^{a-b}+1=3$, $n=2^{a-b}=2$, $q=3\cdot 2^{b-1}-1$ and $h=(q-1)/2$.
   As $h\geq 3$,
   it follows that $q\geq 7$ and $b\geq 3$.
   Therefore, it is routine to check that $q\geq 11$.
   Note that $|G|_2=|\Gamma|_2=2^{b+1}=4(q+1)/3$, we conclude that $G=P\times H$ where $P\in \mathrm{Syl}_{2}(\Gamma)$ and $H=H_0\cong \mathsf{C}_{(q-1)/2}$.
   Thus, $P\cong \mathsf{SD}_{4(q+1)/3}$.
   In fact, as $q=3\cdot 2^{b-1}-1\equiv -1~(\mathrm{mod}~4)$, $P$ is semidihedral by \cite[Page 142]{carterfong1964}.
   So, (5) holds.
 \end{proof}

\begin{hy}\label{hy: rank 3}
  Let $G$ be a solvable primitive permutation group of rank $3$ with $\mathrm{O}_p(G)=1$.
  Write $G=V \rtimes H$, where $V=\mathrm{F}(G)$ is a minimal normal subgroup of $G$ and $H$ is a complement of $V$ in $G$.
  Assume that $|\mathrm{IBr}_{p}^{1}(G)|=2$.
\end{hy}

\begin{lem}\label{lem: rank 3 0}
  Assume Hypothesis \ref{hy: rank 3}.
  Then, for all $v\in V^\sharp$ and $\alpha\in \mathrm{IBr}_{p}(V)^\sharp$,
  $H'$, $\mathrm{C}_{H}(v)$ and $\mathrm{I}_{H}(\alpha)$ are $p$-groups such that $|H:\mathrm{C}_{H}(v)|>1$ and $|H:\mathrm{I}_{H}(\alpha)|>1$.
  Moreover, $|\mathrm{IBr}_{p}^{1}(G|\alpha)|=1$ for all $\alpha \in \mathrm{IBr}_{p}(V)^\sharp$.
\end{lem}
\begin{proof}
  Let $v\in V^\sharp$ and $\alpha \in \mathrm{Irr}(V)^\sharp$.
  As the action of $H$ on $V^\sharp$ has two orbits,
  the action of $H$ on $\mathrm{Irr}(V)^\sharp$ also has two orbits by \cite[Corollary 6.33]{isaacsbook}.
  Since $V$ is a faithful irreducible $H$-module, $|H:\mathrm{C}_{H}(v)|>1$.
  Note that $V$ is a faithful irreducible $H$-module if and only if $\mathrm{Irr}(V)$ is a faithful irreducible $H$-module,
  and so $|H:\mathrm{I}_{H}(\alpha)|>1$.
   Observe that $\mathrm{O}_p(G)=1$, and so $V$ is a $p'$-group. 
  In particular, $\mathrm{IBr}_{p}(V)=\mathrm{Irr}(V)$.
  So, for $\varphi \in \mathrm{IBr}_{p}(G|\alpha)$, $\varphi(1)\ge |H:\mathrm{I}_{H}(\alpha)| >1$ by Clifford's theorem \cite[Corollary 8.7]{navarrobook}.
 As $\mathrm{IBr}_p(V)^\sharp$ is a union of two distinct nontrivial $H$-orbits and $|\mathrm{IBr}_{p}^{1}(G)|=2$,
 it follows that $|\mathrm{IBr}_{p}(G|\alpha)|=|\mathrm{IBr}_{p}^{1}(G|\alpha)|=1$ for each $\alpha \in \mathrm{IBr}_p(V)^\sharp$.
 Also, $|\mathrm{IBr}_{p}^{1}(H)|=|\mathrm{IBr}_{p}^{1}(G/V)|=0$.
 Consequently, $H/\mathrm{O}_{p}(H)$ is an abelian $p'$-group, i.e. $H'\leq \mathrm{O}_{p}(H)$.
 Given that $G=V \rtimes H$, by part (2) of Lemma \ref{lem: linearext}, we deduce that 
 $|\mathrm{IBr}_{p}(\mathrm{I}_{H}(\alpha))|=|\mathrm{IBr}_{p}(G|\alpha)|=1$.
 Therefore, $\mathrm{I}_{H}(\alpha)$ is a $p$-group for each $\alpha \in \mathrm{IBr}_p(V)^\sharp$.
 Furthermore, applying Lemma \ref{lem: stabdual}, we conclude that $\mathrm{C}_{H}(v)$ is a $p$-group for each $v\in V^\sharp$.
\end{proof}

Recall that, for an odd positive integer $n$, the dicyclic group $\mathsf{Dic}_n$ is isomorphic to $\mathsf{C}_{n}  \rtimes \mathsf{C}_{4}$ where $\mathsf{C}_{4}$ acts by inversion on $\mathsf{C}_{n}$.

\begin{prop}\label{prop: rank3 1}
  Assume Hypothesis \ref{hy: rank 3}.
  If $V$ is an $n$-dimensional primitive $H$-module over a prime field $\mathbb{F}_q$, then either $p=2$, $G\cong\mathsf{PrimitiveGroup}(q^n,i)$ where $(q^n,i)\in \{(7^2,18),(17^2,82),(23^{2},51)\}$,
  or $($up to permutation isomorphism$)$ $H\leq \Gamma(V)=\Gamma(q^n)$ and one of the following holds:
  \begin{description}
    \item[(1)] $G$ is a Frobenius group with complement $H$ of order $(q^{n}-1)/2$ such that $H'$ is $p$-group. In this case, one of the following holds:
    \begin{description}
      \item[(1a)] $H\cong \mathsf{C}_{(q^n-1)/2}$;  
      \item[(1b)] $p=n=2$, and $H\cong \mathsf{Q}_{q+1}\times \mathsf{C}_{(q-1)/2}$ where $q \in \mathfrak{M}$ is larger than $3$;
      \item[(1c)] $p$ is a primitive prime divisor of $q^n-1$, and $H'$ is a nontrivial $p$-group;
    \end{description}
    \item[(2)] $p=n=2$, $H\cong \mathsf{D}_{2(q+1)}\times \mathsf{C}_{(q-1)/2}$ 
    where $q\in \mathfrak{M}$;
    \item[(3)] $p=2$, $H\in \mathrm{Syl}_{2}(\Gamma(q^n))$ where $q^n\in \{ 3^{4}, 5^{2} \}$;
    \item[(4)] $p=n=2$, $H\cong \mathsf{SD}_{4(q+1)/3}\times \mathsf{C}_{(q-1)/2}$ where $q=3\cdot 2^{k}-1\geq 11$, and $\mathrm{C}_{H}(v_0)=1$ for some $v_0\in V^\sharp$.
  \end{description}
\end{prop}
\begin{proof}
  Note that $G$ is a solvable primitive permutation group of rank $3$ and that $V$ is a primitive $H$-module.
  Applying Theorem \ref{thm: classification of solvable primitive permutation groups of rank 3} and checking 
  via $\mathsf{GAP}$ \cite{gap},
  we have that either $p=2$, $G\cong\mathsf{PrimitiveGroup}(q^n,i)$ where $(q^n,i)\in \{(7^2,18),(17^2,82),(23^{2},51)\}$,
  or (up to permutation isomorphism) $H\leq \Gamma(q^n)$.
  If $H$ is abelian, as $V$ is a faithful irreducible $H$-module, $H\leq \Gamma_0(q^{n})$.
  So, (1a) holds.
  If $H$ is nonabelian, then our desired results follow directly by Lemmas \ref{lem: rank 3 0} and \ref{lem: rank3hgamma}.
\end{proof}

\begin{rmk}
  {\rm There are groups satisfying  (1c) of Proposition \ref{prop: rank3 1}.
  For instance, 
  let $G$ be a Frobenius group with kernel $V\cong (\mathsf{C}_{11})^{2}$ and 
  nonabelian complement $H\cong \mathsf{C}_{5}\times \mathsf{Dic}_3$. 
  The group $G$ is known as $\mathsf{PrimitiveGroup}(11^2,31)$.
  It is routine to verify that $3$ is a primitive prime divisor of $11^2-1$, $H\leq \Gamma(11^2)$, $H'\cong \mathsf{C}_{3}$ and $|\mathrm{IBr}_{3}^{1}(G)|=2$.}
\end{rmk}

\begin{prop}\label{prop: rank3 2}
  Assume Hypothesis \ref{hy: rank 3}.
  If $V$ is an $n$-dimensional imprimitive $H$-module over a prime field $\mathbb{F}_q$, then $V=V_1\times V_2$ with $|V_i|=q^k$ and,
  $($up to permutation isomorphism$)$
  $H\leq T\wr \langle \sigma\rangle$ where $T\leq \Gamma(V_1)=\Gamma(q^{k})$ and $\sigma$ is an involution transposing $V_1$ and $V_2$, and we are in one of the following cases:
  \begin{description}
    \item[(1)] $p=2$, $q=q^k \in \mathfrak{F}$, $T=\Gamma_0(q)=\langle t\rangle$ and, either $H= T \wr \langle \sigma\rangle$ or $H=\langle t \sigma\rangle \Phi(T\wr \langle \sigma\rangle) $;
    \item[(2)] $p=2$, $G\cong\mathsf{PrimitiveGroup}(3^4,i)$ where $i\in  \{36,37, 41, 61,62,63,65,66,67, 82,84,85,86,87,95\}$;
    \item[(3)] $p>2$, $T=\Gamma_0(q^{k})$ and $H\cong \mathsf{Dic}_{(q^{k}-1)/2}\times \mathsf{C}_{(q^{k}-1)/2}$ where $q^k=2p^s+1$ such that $$(q^{k},p,s)\in \{ (q,p,s), (3^5,11,2), (3^k,(3^k-1)/2,1) \}.$$
  \end{description}
\end{prop}
\begin{proof}
 By Lemma \ref{lem: rank 3 0}, we know that $H'$ and $\mathrm{C}_{H}(v)$ are $p$-groups for all $v\in V^\sharp$.
 In particular, $H/\mathrm{O}_{p}(H)$ is an abelian $p'$-group.
 Since $V$ is a $p'$-group, $\mathrm{IBr}_{p}(V)=\mathrm{Irr}(V)$.
 Note that $G$ is a solvable primitive permutation group of rank $3$ and that $V$ is an imprimitive $H$-module,
 and so $H$ is nonabelian.
  Applying Theorem \ref{thm: classification of solvable primitive permutation groups of rank 3},
  we have that $V=V_1\times V_2$, $\mathrm{N}_H(V_i)$ acts transitively on $V_i^\sharp$, and $V_1^\sharp\cup V_2^\sharp$ and $V-(V_1\cup V_2)$ are the nontrivial orbits of the action of $H$ on $V$.
 So, \cite[Lemma 2.8]{manzwolfbook} yields that $H$ is isomorphic to a subgroup of $T \wr \langle \sigma\rangle$ where $T=\mathrm{N}_{H}(V_1)/\mathrm{C}_{H}(V_1)$ and $\sigma$ is an involution transposing $V_1$ and $V_2$.
  In the following, we identify $H$ as a subgroup of $T \wr \langle \sigma\rangle$.
  As $\mathrm{N}_{H}(V_i)$ has index $2$ in $H$, we have that $\mathrm{N}_{H}(V_1)=\mathrm{N}_{H}(V_2)$ is normal in $H$.
  Set $N=\mathrm{N}_{H}(V_i)$ and $Y_i=V_i \rtimes T_i$ where $T_i=N/\mathrm{C}_{H}(V_i)$ for $i=1,2$.
  Then $Y_i$ are solvable doubly transitive permutation group with $\mathrm{O}_{p}(Y_i)=1$.
  Next, we proceed the rest proof by verifying seven claims.

  \textbf{Claim 1.} $|V_i|>2$ and $T_i>1$.

  Otherwise, $|H|\leq 2$ whereas $H$ is nonabelian.

  \textbf{Claim 2.} $T_i\leq \Gamma(V_i)=\Gamma(q^{k})$. 
  In this case, either $Y_i$ is a doubly transitive Frobenius group with complement $T_i$ such that $T_i'$ is a $p$-group, or $p=k=2$, $Y_i=A\Gamma(q^{2})$ where $q \in \mathfrak{M}$.
  
  Note that $H/\mathrm{O}_{p}(H)$ is an abelian $p'$-group and that $T_i=N/\mathrm{C}_{H}(V_i)$ where $N\unlhd H$,
  and so $T_i/\mathrm{O}_{p}(T_i)$ is also an abelian $p'$-group.
  This implies that $|\mathrm{IBr}_{p}^{1}(T_i)|=0$.
  Let $v_i \in V_i^\sharp$, and observe that 
   $\mathrm{C}_{T_i}(v_i)=\mathrm{C}_{N/\mathrm{C}_{H}(V_i)}(v_i)=\mathrm{C}_{N}(v_i)/\mathrm{C}_{H}(V_i)\leq \mathrm{C}_{H}(v_i)/\mathrm{C}_{H}(V_i)$ and that $\mathrm{C}_{H}(v_i)$ is a $p$-group,
   and thus $\mathrm{C}_{T_i}(v_i)$ is a $p$-group.
   Let $\alpha_i\in \mathrm{IBr}_p(V_i)^\sharp$.
   Since $V_i$ is a $p'$-group, $\mathrm{IBr}_{p}(V_i)=\mathrm{Irr}(V_i)$.
   By Lemma \ref{lem: stabdual}, we know that $\mathrm{I}_{T_i}(\alpha_i)$ is also a $p$-group.
   As $T_i$ acts transitively on $V_i^\sharp$, it also acts transitively on $\mathrm{IBr}_p(V_i)^\sharp$.
   Noting that $|\mathrm{IBr}_p(V_i)^\sharp|=|V_i^\sharp|>1$ by Claim 1 and applying Clifford's theorem \cite[Corollary 8.7]{navarrobook}, we deduce that $\mathrm{IBr}_{p}(Y_i|\alpha_i)=\mathrm{IBr}_{p}^{1}(Y_i)$.
   Since $Y_i=V_i  \rtimes T_i$ such that $\mathrm{I}_{T_i}(\alpha_i)$ is a $p$-group,
   we deduce by part (2) of Lemma \ref{lem: linearext} that $|\mathrm{IBr}_{p}(Y_i|\alpha_i)|=|\mathrm{IBr}_{p}(\mathrm{I}_{T_i}(\alpha_i))|=1$.
   Consequently, $|\mathrm{IBr}_{p}^{1}(Y_i)|=1$ and $\mathrm{O}_{p}(Y_i)=1$.

   We now assume that $p=2$ and $Y_i\cong\mathsf{PrimitiveGroup}(q^k,j)$, where $(q^{k},j)\in \{ (5^2,18),(3^4,71) \}$.
   Let $a$ and $b$ be the sizes of the two nontrivial $H$-orbits in $V$ such that $a\leq b$.
   Then $a=|V_1^\sharp\cup V_2^\sharp|=2(q^k-1)$ and $b=|V|-|V_1\cup V_2|=(q^k-1)^{2}$.
   Recall that $\mathrm{C}_{H}(v)$ is a $2$-group for each $v\in V^\sharp$,
   and hence it follows that $a/a_2=b/b_2$.
   However, this will not happen because $q^k\in \{ 5^2, 3^4 \}$.
   So, by \cite[Theorem A]{dolfi}, Claim 2 holds.

   \textbf{Claim 3.} Let $x\in H-N$, $\lambda_1\in \mathrm{Irr}(V_1)^\sharp$ and $\lambda_2=\lambda_1^x \in \mathrm{Irr}(V_2)^\sharp$.
   Set $\lambda=\lambda_1\times \lambda_2$ and $\mu=\lambda_1\times 1_{V_2}$.
   Then $\lambda$ and $\mu$ lie in distinct $H$-orbits.
   Assume that $\mathrm{I}_{T_i}(\lambda_i)=1$ for $1\leq i\leq 2$. 
  Then $\mathrm{I}_{N}(\lambda_i)=\mathrm{C}_{H}(V_i)$, $\mathrm{I}_{H}(\mu)=\mathrm{C}_{H}(V_1)$ and $|\mathrm{I}_{H}(\lambda)|\leq 2$. 
  If in addition $\mathrm{I}_{H}(\lambda)=1$, then there is no involution in $H-N$.

   It is routine to check that $\mathrm{I}_{N}(\lambda)=\mathrm{I}_{N}(\lambda_1)\cap \mathrm{I}_{N}(\lambda_2)$ and $\mathrm{I}_{H}(\mu)=\mathrm{I}_{N}(\lambda_1)$, and that $\lambda$ and $\mu$ are not $H$-conjugate.
   Also, as $|H:N|=2$, we have $|\mathrm{I}_{H}(\lambda):\mathrm{I}_{N}(\lambda)|\leq 2$.
   Assume that $\mathrm{I}_{T_i}(\lambda_i)=1$ for $1\leq i\leq 2$.
   Since $\mathrm{I}_{T_i}(\lambda_i)=\mathrm{I}_{N}(\lambda_i)/\mathrm{C}_{H}(V_i)$,
   $\mathrm{I}_{N}(\lambda_i)=\mathrm{C}_{H}(V_i)$.
   So, $\mathrm{I}_{H}(\mu)=\mathrm{C}_{H}(V_1)$ and $\mathrm{I}_{N}(\lambda)=\mathrm{C}_{H}(V_1)\cap \mathrm{C}_{H}(V_2)=1$.
   Observe that $|\mathrm{I}_{H}(\lambda):\mathrm{I}_{N}(\lambda)|\mid 2$,
   and hence $|\mathrm{I}_{H}(\lambda)|\leq 2$.
   Note that all involutions in $H-N$ lie in $\mathrm{I}_{H}(\lambda)$,
   and so the last statement of Claim 3 holds.

 \textbf{Claim 4.} 
  Let $\mathcal{O}_1=\{ \alpha\times \beta: \alpha \in \mathrm{Irr}(V_1)^\sharp,~\beta \in \mathrm{Irr}(V_2)^\sharp \}$ and $\mathcal{O}_2=\mathrm{Irr}(V)^\sharp-\mathcal{O}_1$.
  Then $\mathcal{O}_1$ and $\mathcal{O}_2$ are the orbits of the action of $H$ on $\mathrm{Irr}(V)^\sharp$ containing $\lambda$ and $\mu$ respectively, such that $|\mathcal{O}_1|=(q^k-1)^{2}$ and $|\mathcal{O}_2|=2(q^{k}-1)$.
  Moreover, $q^k-1=2p^s$ where $q>2$ and $p^s=|\mathrm{I}_{H}(\mu)|/|\mathrm{I}_{H}(\lambda)|$, and $p$ is a primitive prime divisor of $q^{k}-1$ if $p>2$.

  As the action of $H$ on $V^\sharp$ has two orbits,
  the action of $H$ on $\mathrm{Irr}(V)^\sharp$ also has two orbits by \cite[Corollary 6.33]{isaacsbook}.
  So, by Claim 3, $\mathrm{Irr}(V)^\sharp=\mathcal{O}_\lambda\sqcup \mathcal{O}_\mu$ where $\mathcal{O}_\lambda=\{ \lambda^h:h\in H \}$ and $\mathcal{O}_\mu=\{ \mu^h : h\in H \}$.
  Since $N$ acts transitively on both $V_1^\sharp$ and $V_2^\sharp$, it also acts transitively on both $\mathrm{Irr}(V_1)^\sharp$ and $\mathrm{Irr}(V_2)^\sharp$.
  Therefore, it is routine to check that $\mathcal{O}_\lambda=\mathcal{O}_1$, $\mathcal{O}_\mu=\mathcal{O}_2$ and $|\mathcal{O}_1|=(q^k-1)^{2}$ and $|\mathcal{O}_2|=2(q^{k}-1)$.
  Recall that $\mathrm{C}_{H}(v)$ are $p$-groups for all $v\in V^\sharp$, 
  and so $\mathrm{I}_{H}(\alpha)$ are $p$-groups for all $\alpha \in \mathrm{Irr}(V)^\sharp$ by Lemma \ref{lem: stabdual}. 
  In particular, $|\mathrm{I}_{H}(\mu)|/|\mathrm{I}_{H}(\lambda)|=p^s$ for some nonnegative integer $s$.
  Since $|H:\mathrm{I}_{H}(\lambda)|=|\mathcal{O}_1|=(q^{k}-1)^{2}$ and $|H:\mathrm{I}_{H}(\mu)|=|\mathcal{O}_2|=2(q^{k}-1)$,
  we conclude that $q^k-1=2p^s$.
  If $p>2$, then, according to Lemma \ref{lem: zsig}, $q^k-1$ has a primitive prime divisor $p$.

\textbf{Claim 5.} $4p^{2s}\mid |H|\mid 8k^2 p^{2s}$.

Since $H\leq T\wr \langle \sigma\rangle$ where $o(\sigma)=2$ and $T\leq \Gamma(q^{k})$, 
$|H|\mid 2|\Gamma(q^{k})|^2$.
Note that $q^k-1=2p^s$, and so $|H|\mid 8k^2 p^{2s}$.
Also, $4p^{2s}=|\mathcal{O}_1|\mid |H|$ by Claim 4.

\textbf{Claim 6.} If $p=2$, then either (1) or (2) holds.

Assume that
$p=2$.
Note that $q^{k}=2^{s+1}+1$,
and so Lemma \ref{lem: weakcatalan} implies that either $q=q^k=2^{s+1}+1\in \mathfrak{F}$ or $q^k=3^2$ and $2^{s+1}=2^3$.
Also, by Claim 5, we have $2^{2s+2}\mid |H|\mid 2^{2s+3}k^{2}$.

If $q=q^k=2^{s+1}+1$ is a Fermat prime, then $Y_i\cong A\Gamma(V_i)=A\Gamma(q)$ is a doubly transitive Frobenius group with complement $T_i=\Gamma_0(q)\cong \mathsf{C}_{q-1}$.
In particular, $\mathrm{I}_{T_i}(\lambda_i)=1$.
By Claim 3, we have that $\mathrm{I}_{H}(\mu)=\mathrm{C}_{H}(V_1)$ and $|\mathrm{I}_{H}(\lambda)|\le 2$.
 If $|\mathrm{I}_{H}(\lambda)|=2$, then $|H|=2(q-1)^2$.
 As $|T \wr \langle \sigma\rangle|=2(q-1)^2$,
 we conclude that $H= T \wr \langle \sigma\rangle$.
 If $\mathrm{I}_{H}(\lambda)=1$, then $|H|=(q-1)^2$ and there is no involution in $H-N$ by Claim 3.
 As $|T \wr \langle \sigma\rangle|=2(q-1)^2$, $H$ is a maximal subgroup of the 2-group $T \wr \langle \sigma\rangle$.
 Moreover, since $|H|=2(q-1)|\mathrm{I}_{H}(\mu)|$ and $\mathrm{I}_{H}(\mu)=\mathrm{C}_{H}(V_1)$, we deduce that $|\mathrm{C}_{H}(V_1)|=(q-1)/2$.
 Recall that $H$ is nonabelian, and so $|H|=(q-1)^{2}\geq 8$.
 Therefore, $2\mid (q-1)/2$.
  Given that $N/\mathrm{C}_{H}(V_i)=T_i\cong \mathsf{C}_{q-1}$ and that $\mathrm{C}_{H}(V_1)\cap \mathrm{C}_{H}(V_2)=1$,
we have $[N,\mathrm{C}_{H}(V_1)\mathrm{C}_{H}(V_2)]\leq \mathrm{C}_{H}(V_1)\cap \mathrm{C}_{H}(V_2)=1$.
Since $|N:\mathrm{C}_{H}(V_1)\mathrm{C}_{H}(V_2)|= 2$ and $\mathrm{C}_{H}(V_1)\mathrm{C}_{H}(V_2)\cong (\mathsf{C}_{(q-1)/2})^{2}$ is central in $N$,
$N$ is an abelian $2$-group isomorphic to 
$ \mathsf{C}_{q-1}\times \mathsf{C}_{(q-1)/2}$.
In particular, $N$ has only 3 involutions.
So, $H$ has exactly 3 involutions.
Recall that $H$ is a nonabelian maximal subgroup of $T \wr \langle \sigma\rangle$,
and hence it follows by Lemma \ref{lem: wreath2} that 
$H=\langle t \sigma\rangle \Phi(T \wr \langle \sigma\rangle)$ where $T=\Gamma_0(q)=\langle t\rangle$.
So, (1) holds.

If $q^k=3^2$ and $2^{s+1}=2^3$, then $2^{6}\mid |H|\mid 2^{9}$.
Note that $G$ is a solvable primitive permutation group of degree $3^4$ such that
$|\mathrm{IBr}_{2}^{1}(G)|=2$ and $2^6\leq |G/\mathrm{F}(G)|=|H|\leq 2^9$.
Checking the groups with this property in $\mathsf{GAP}$ \cite{gap},
we conclude that (2) holds.

\textbf{Claim 7.} If $p>2$, then (3) holds.

 Since $p>2$, we deduce by Claim 2 that $T_i\leq \Gamma(q^k)$ and $Y_i$ is a doubly transitive Frobenius group 
 with complement $T_i$ such that $T_i'$ is a $p$-group. 
 Furthermore, as the Frobenius complement $T_i$
 has order $q^k-1=2p^s$ where $p>2$, 
we have $T_i=\Gamma_0(q^{k})\cong \mathsf{C}_{q^k-1}$.
In fact, as $p$ is a primitive prime divisor of $q^{k}-1$ by Claim 4, \cite[Corollary 6.6, Lemma 6.4]{manzwolfbook} forces that $P_i\leq \Gamma_0(q^{k})$ for $P_i\in \mathrm{Syl}_{p}(T_i)$;
note that $q^k-1=2p^s$ where $p>2$, and so $k$ is odd which implies that $R_i\leq \Gamma_0(q^{k})$ for $R_i\in \mathrm{Syl}_{2}(T_i)$; hence $T_i=\Gamma_0(q^{k})\cong \mathsf{C}_{q^k-1}$.

Since $T_i$ acts fixed-point-freely on $V_i$, $\mathrm{I}_{T_i}(\lambda_i)=1$.
So, Claim 3 shows that $\mathrm{I}_{H}(\mu)=\mathrm{C}_{H}(V_1)$ and $|\mathrm{I}_{H}(\lambda)|\leq 2$.
As $\mathrm{I}_{H}(\lambda)$ is a $p$-group where $p>2$, we deduce that $\mathrm{I}_{H}(\lambda)=1$.
Consequently, again by Claim 3, $H-N$ contains no involution.
Also, as $|H|=|\mathcal{O}_1||\mathrm{I}_{H}(\lambda)|=|\mathcal{O}_2||\mathrm{I}_{H}(\mu)|$,
it follows by Claim 4 that $|H|=(q^k-1)^2=2(q^k-1)|\mathrm{C}_{H}(V_1)|$, and so $|\mathrm{C}_{H}(V_i)|=(q^k-1)/2$.
Since $H\leq T \wr \langle \sigma\rangle$ and $|T \wr \langle \sigma\rangle|=2(q^k-1)^2$,
$H$ has index 2 in $T \wr \langle \sigma\rangle$.

Set $C_i=  \mathrm{C}_{H}(V_i)$ and $C=C_1\times C_2$.
Since $(|H/C|,|C|)=1$ and $|C|=p^{2s}$,
we conclude, by the Schur-Zassenhaus theorem, that $H=C  \rtimes R$ where $C\in \mathrm{Syl}_{p}(H)$ and $R\in \mathrm{Syl}_{2}(H)$.
Note that $N=C  \rtimes (R\cap N)$ and that
$H-N$ has no involution, and so $R\cong \mathsf{C}_{4}$.
Next, we claim that $H\cong \mathsf{Dic}_{(q^k-1)/2}\times \mathsf{C}_{(q^k-1)/2}$.
Write $C_i=\langle c_i\rangle$ and $R=\langle y\rangle$,
since $N=C \rtimes \langle y^2\rangle$ and $N/C_i\cong \mathsf{C}_{q^k-1}$,
we have 
$[y^2,N]\leq C_1\cap C_2=1$,
i.e. $N=C \times \langle y^2\rangle$.
Also, direct calculations yield that $C\cap \mathrm{Z}(H)=\langle c_1c_1^y\rangle\cong \mathsf{C}_{(q^k-1)/2}$ and $H'=\langle [c_1,y]\rangle\cong \mathsf{C}_{(q^k-1)/2}$.
Since $y$ acts by inversion on $H'$, $H=(C\cap \mathrm{Z}(H))\times (H' \rtimes R)$ where $C\cap\mathrm{Z}(H)\cong\mathsf{C}_{(q^k-1)/2}$ and $H'R\cong \mathsf{Dic}_{(q^k-1)/2}$.

Finally, we show that $(q^{k},p,s)\in \{ (q=2p^s+1,p,s), (3^5,11,2), (3^{k},(3^k-1)/2,s)\}$.
To establish this, we may assume that $k>1$.
Recall that $q^k-1=2p^s$ where $p>2$,
and so $k$ is odd and, by Claim 4, $p$ is a primitive prime divisor of $q^k-1$.
In particular, $k<p$.
Since $(q-1)(q^{k-1}+\cdots +q+1)=2p^{s}$ and $p\nmid q-1$,
we conclude that $q-1=2$ and $q^{k-1}+\cdots +q+1=p^{s}$.
So, by calculation, $q=3$ and $3^k-2p^s=1$.
Observe that $k\geq 3$ and $p>2$,
and hence it follows by \cite[Corollary 1]{bm99} that either $(q^k,p,s)=(3^5,11,2)$ or $s=1$ and $p=(3^k-1)/2$.
\end{proof}

\subsection{Solvable groups with nonabelian Fitting subgroup}

In this subsection, we deal with solvable groups with nonabelian Fitting subgroups that have 
exactly two nonlinear irreducible $p$-Brauer characters.

\begin{lem}\label{lem: structure}
  Let a solvable group $H$ act via automorphisms on a nonabelian $q$-group $Q$, and let $G=Q \rtimes H$.
  Set $\overline{G}=G/\mathrm{Z}(Q)$.
  Assume that $\overline{G}$ is a doubly transitive Frobenius group with complement $\overline{H}$ such that $\overline{H}'$ is a $p$-group, and that $H$ acts transitively on $\mathrm{Z}(Q)^\sharp$.
  Set $|\overline{Q}|=q^{n}$ and $|\mathrm{Z}(Q)|=q^{d}$.
  If $(q^{n}-1)=p^k(q^{d}-1)$ and $d<n$,
  then $n=2d$, $|H|=q^{n}-1$
  and one of the following holds.
  \begin{description}
    \item[(1)] $q=2$, $p=3$ and $G\cong \mathrm{SL}_2(3)$.
    \item[(2)] $q=2$, $p\mid 2^d+1$ and $G$ is isomorphic to $B(2^{d})$ where $2^{d}+1~(>3)$ is either a Fermat prime or $9$.
    \item[(3)] $q\in \mathfrak{M}$, $p=n=2$,  
    $Q$ is isomorphic to $\mathsf{ES}(q^{3}_+)$, and $H$ is a maximal subgroup of $\Gamma(q^{2})$ which 
    is isomorphic to either $\mathsf{C}_{q^{2}-1}$, or $\mathsf{Q}_{2(q+1)}\times \mathsf{C}_{(q-1)/2}$ with $q\geq 7$.
  \end{description}
\end{lem}
\begin{proof}
    Set $Z=\mathrm{Z}(Q)$.
    Observe that $\overline{G}$ is a doubly transitive Frobenius group with complement $\overline{H}$,
    and that $H$ acts transitively on $Z^\sharp$.
    So, we deduce that $|H|=|\overline{H}|=q^{n}-1$,
    and the nonabelian $q$-group $Q$ is special.
    In particular, $(|H|,|Q|)=1$.
    Further, Lemma \ref{lem: omega} implies that $Z=\Omega_1(Q)$ if $q=2$.
    Note that $q^{d}-1\mid q^{n}-1$ where $d<n$, and so it follows that $d$ is a proper divisor of $n$.

    Let $P$ be a Sylow $p$-subgroup of $H$.
    As $\overline{H}'\leq \overline{P}$, we have $P\unlhd H$.
    Hence, the Schur-Zassenhaus theorem implies that $H=P \rtimes L$.
    Since $H$ is isomorphic to a Frobenius complement $\overline{H}$,
    $P$ is either generalized quaternion or cyclic, and $L$ is cyclic.

    Assume that $p$ is a primitive prime divisor of $q^{n}-1$.
    As $d<n$, $(p,q^d-1)=1$ and $p>2$.
    By \cite[Proposition 6.3]{manzwolfbook}, we know that $P$ is cyclic.
    Note that $H$ acts transitively on $Z^\sharp$ and $q^n-1=p^k(q^d-1)$, and that $P\unlhd H$.
    So, it follows that $P=\mathrm{C}_{H}(z)=\mathrm{C}_{H}(Z)\cong \mathsf{C}_{p^{k}}$ for each $z\in Z^\sharp$.
    Now, let $C=\mathrm{C}_{G}(Z)$.
    Then $C=Q \rtimes P$ and $\mathrm{Z}(C)=Z$.
    Let $K$ be a maximal subgroup of $Z$ and set $\widetilde{C}=C/K$.
    Working in the group $\widetilde{C}$,
    we see that $\widetilde{C}=\widetilde{Q} \rtimes \widetilde{P}$.
    As $\overline{Q}\cong (\mathsf{C}_{q})^{n}$ is an irreducible $H$-module and $p$ is a primitive prime divisor of $q^{n}-1$,
    $\overline{Q}$ is $P$-irreducible.
    Note that $|\widetilde{Z}|=q$,
    and so $\widetilde{Q}/\widetilde{Z}$ and $\widetilde{Z}$ are $\widetilde{C}$-chief factors in the nonabelian $q$-group $\widetilde{Q}$.
    Therefore, $\widetilde{Q}$ is an extraspecial $q$-group.
    Applying \cite[Corollary 2]{winter1972} to $\widetilde{P}$ and $\widetilde{Q}$, we conclude that $|\widetilde{P}|\mid q^{n/2}+1$.
    Since $|\widetilde{P}|=|P|=(q^{n}-1)/(q^{d}-1)$ where $d$ is a proper divisor of $n$,
    it is routine to check that $n=2d$.
    Further, $p^{k}=q^{n/2}+1$.
    As $p$ is an odd primitive prime divisor of $q^{n}-1$,
    by Lemma \ref{lem: weakcatalan}, $q=2$ and $p=p^{k}=2^{n/2}+1$ is a Fermat prime. 
    Observe that $L/\mathrm{C}_{L}(P)$ is isomorphic to a $2'$-subgroup of the $2$-group $\mathrm{Aut}(P)\cong \mathrm{Aut}(\mathsf{C}_{p})$.
     Therefore, $L=\mathrm{C}_{L}(P)$, and thus $H=P\times L$ is a cyclic group of order $2^{n}-1$.
    Recall that as $q=2$, $Z=\Omega_1(Q)$. 
    If $n=2$, then $Q\cong \mathsf{Q}_8$,
    and so $p=3$ and $G\cong \mathrm{SL}_2(3)$.
    Suppose that $n>2$.
    As now $H\cong \mathsf{C}_{2^{n}-1}$, who is isomorphic to a subgroup of $\mathrm{Aut}(Q)$, acts transitively on $\Omega_1(Q)^\sharp$ and $|Q|=|Z|^3=2^{3d}$ where $d=n/2$ is a power of $2$,
    it follows by Lemma \ref{lem: suzuki2} that $Q$ is a Suzuki $2$-group of $B$-type.
    Consequently, this leads us to conclude, by Proposition \ref{prop: suzuki2}, that $G$ is isomorphic to $B(2^{d})$ where $2^{d}+1\in \mathfrak{F}$ and $d>1$.

    Assume now that $p$ is not a primitive prime divisor of $q^{n}-1$.
    Since $q^{n}-1=p^{k}(q^{d}-1)$ and $ d<n$, $q^{n}-1$ does not have a primitive prime divisor.
    An application of Lemma \ref{lem: zsig} to $q^{n}-1$ yields that either $n=2$ and $q \in \mathfrak{M}$, or $q^{n}=2^{6}$.
   
    If $q^{n}=2^{6}$, then $d=3$, $p^{k}=3^{2}$
    and $|H|=q^{n}-1=63$.
    In this case, $Q$ is a special $2$-group with $Z=\Omega_1(Q)$ and $P\cong \mathsf{C}_{9}$.
    Observe that $L/\mathrm{C}_{L}(P)$ is isomorphic to a $\{ 2,3 \}'$-subgroup of $\mathrm{Aut}(P)\cong \mathsf{C}_{6}$,
    and hence $L=\mathrm{C}_{L}(P)$ which implies $H=P\times L\cong \mathsf{C}_{63}$.
    As $H$, who is isomorphic to a subgroup of $\mathrm{Aut}(Q)$, acts transitively on $\Omega_1(Q)^\sharp$ and $|Q|=|Z|^3=2^{9}$,
    $Q$ is a Suzuki $2$-group of either $B$-type or $C$-type (see \cite[Page 82]{higman1963}).
    Since $H=P\times L$ where $L\cong \mathsf{C}_{7}$,
    Clifford's theorem \cite[Theorem 0.1]{manzwolfbook} implies that $\overline{Q}$ is a homogeneous $L$-module.
    Noting that 
    $L$ acts transitively on $\Omega_1(Q)^\sharp$, 
    we deduce
    by Lemma \ref{lem: suzuki2}
    that $Q$ is a Suzuki $2$-group of $B$-type.
    As a consequence, this leads us to conclude, by Proposition \ref{prop: suzuki2}, that $G$ is isomorphic to $B(8)$.

    Assume that $n=2$ and $q$ is a Mersenne prime.
    Then $d=1$ and $Q$ is an extraspecial $q$-group of order $q^{3}$.
    Moreover, $p^{k}=q+1=2^{k}$ where $k$ is a prime.
    Since $H$ acts transitively on both $\overline{Q}^\sharp$ and $Z^\sharp$,
    we deduce by Lemma \ref{lem: omega} that $Q\cong \mathsf{ES}(q^{3}_+)$.
   Recall that $P$ is either a generalized quaternion $2$-group or a cyclic $2$-group, and that $|P|=(q^{2}-1)_2=2(q+1)=2^{k+1}$.
     It is also worth noting that $\mathrm{Aut}(P)$ is a $2$-group when $P$ is not isomorphic to $\mathsf{Q}_8$, and that $\mathrm{Aut}(\mathsf{Q}_8)\cong \mathsf{S}_4$.
    As $H$ acts transitively on $Z^\sharp$, we conclude that $P$ is not isomorphic to $\mathsf{Q}_8$.
    Indeed, if $P\cong \mathsf{Q}_8$, then $2^{k+1}=2(q+1)=8$,
    so direct calculations imply that $Q\cong\mathsf{ES}(3^3_+)$ and $H\cong\mathsf{Q}_8$;
    however, in this case, $[H,Z]=1$, a contradiction.  
    Therefore, $\mathrm{Aut}(P)$ is a $2$-group and $P$ is not isomorphic to $\mathsf{Q}_8$.
    Observing that $L/\mathrm{C}_{L}(P)$ is isomorphic to a $2'$-subgroup of the $2$-group $\mathrm{Aut}(P)$,
    we conclude that $H=P\times L$ where $|P|=2(q+1)$ and $L\cong \mathsf{C}_{(q-1)/2}$.
    Now, let $H_1$ be a cyclic subgroup of $H$ of order $(q^2-1)/2$. 
    Noting that $H_1$ acts faithfully and irreducibly on $\overline{Q}$, by \cite[Theorem 2.1]{manzwolfbook}, we deduce that $H$ is isomorphic to a subgroup of $\Gamma(\overline{Q})=\Gamma(q^2)$.
    Now, we identify $H$ as a subgroup of $\Gamma(q^2)$.
     As $q$ is a Mersenne prime, $q\equiv -1~(\mathrm{mod}~4)$,
    and so Sylow $2$-subgroups of $\mathrm{GL}_2(q)$ are semidihedral of order $4(q+1)$
    (see, for instance, \cite[Page 142]{carterfong1964}).
    It is noteworthy that $\Gamma(q^2)\leq \mathrm{GL}_2(q)$ and that $|\Gamma(q^{2})|_2=|\mathrm{GL}_2(q)|_2$,
  we deduce by \cite[Lemma 6.5]{manzwolfbook} that $\Gamma(q^2)=R\times S$ where $R\cong \mathsf{SD}_{4(q+1)}$ and $S\cong \mathsf{C}_{(q-1)/2}$ is of odd order.
   Recalling that $|H|=q^2-1$ and $|\Gamma(q^2)|=2(q^2-1)$,
   we infer that $P$ has index $2$ in $R$ and $L=S$.
   Notably, as the only three maximal subgroups of $\mathsf{SD}_{4(q+1)}$ are $\mathsf{D}_{2(q+1)}$, $\mathsf{C}_{2(q+1)}$ and $\mathsf{Q}_{2(q+1)}$,
    we conclude that $H$ is the maximal subgroup of $\Gamma(q^{2})$ which is isomorphic to either $\mathsf{C}_{q^2-1}$ or $\mathsf{Q}_{2(q+1)}\times \mathsf{C}_{(q-1)/2}$.
    Recall that $P$ is not isomorphic to $\mathsf{Q}_8$,
    and so, if $H\cong \mathsf{Q}_{2(q+1)}\times \mathsf{C}_{(q-1)/2}$, then $q\geq 7$.
  \end{proof}

\begin{prop}\label{fnonabel}
  Let $G$ be a finite solvable group with $\mathrm{O}_p(G)=1$ and set $\Phi=\Phi(G)$.
  Assume that $|\mathrm{IBr}_{p}^{1}(G/\Phi)|=1$.
  If $|\mathrm{IBr}_{p}^{1}(G)|=2$, then one of the following holds:
  \begin{description}
    \item[(1)] $p=3$, $G$ is isomorphic to $\mathrm{SL}_2(3)$;
    \item[(2)] $p\mid 2^{m}+1$, $G$ is isomorphic to $B(2^{m})$
    where either $m=3$ or $p=2^{m}+1>3$ is a Fermat prime;
    \item[(3)] $p=n=2$, $q \in \mathfrak{M}$ and
    $G$ is isomorphic to $\mathsf{ES}(q^{3}_+) \rtimes X$ where either $X=\Gamma(q^{2})$ or, $X$
    is a maximal subgroup of $\Gamma(q^{2})$ which is isomorphic to either  $\mathsf{C}_{q^{2}-1}$, or $\mathsf{Q}_{2(q+1)}\times \mathsf{C}_{(q-1)/2}$ with $q\geq 7$.
  \end{description}
\end{prop}
\begin{proof}
  Set $F=\mathrm{F}(G)$ and $\overline{G}=G/\Phi$.
  Since $\Phi=\Phi(G)$, it follows that $\mathrm{F}(\overline{G})=\overline{F}$.
  As $\mathrm{O}_{p}(G)=1$, $F$ is a $p'$-group.
  In particular, $\mathrm{O}_{p}(\overline{G})=1$, $\mathrm{IBr}_{p}(F)=\mathrm{Irr}(F)$, $\mathrm{IBr}_{p}(\overline{F})=\mathrm{Irr}(\overline{F})$ and $\mathrm{IBr}_{p}(\Phi)=\mathrm{Irr}(\Phi)$.
  Noting that $|\mathrm{IBr}_{p}^{1}(G)|=2$ and $|\mathrm{IBr}_{p}^{1}(\overline{G})|=1$,
  we write $\mathrm{IBr}_{p}^{1}(G)-\mathrm{IBr}_{p}^{1}(\overline{G})=\{ \varphi \}$ and $\mathrm{IBr}_{p}^{1}(\overline{G})=\{ \psi \}$.
  Now, we proceed the proof by proving several claims.

  \smallskip

  \textbf{Claim 1.} $\mathrm{IBr}_{p}(\Phi)^\sharp=\mathrm{IBr}_{p}(\varphi_\Phi)$.
  Moreover, $G$ acts transitively on both $ \mathrm{Irr}(\Phi)^\sharp$ and $\Phi^\sharp$. In particular, $\Phi$ is minimal normal in $G$.

  This claim follows directly by part (2) of Lemma \ref{lem: ibrG/N}.

   \smallskip

\textbf{Claim 2.} $\overline{G}=\overline{F} \rtimes \overline{H}$ is a doubly transitive permutation group with a minimal normal subgroup $\overline{F}=\mathrm{F}(\overline{G})$ such that $|\overline{F}|>2$ and $\mathrm{C}_{\overline{G}}(x)/\overline{F}$ is a $p$-group for each $x \in \overline{F}^\sharp$.
Moreover, one of the following holds: $\overline{G}$ is a doubly transitive Frobenius group with complement $\overline{H}$ such that
$\overline{H}'$ is a $p$-group; $\overline{G}\cong A\Gamma(q^{2})$ where $q \in \mathfrak{M}$ and $p=2$; $\overline{G}\cong\mathsf{PrimitiveGroup}(q^n,i)$ where $(q^{n},i)\in \{ (5^{2},18),(3^{4},71) \}$ and $p=2$.
Also, $F$ is a $G$-indecomposable subgroup.
In particular, $G=F\rtimes H$ where $F$ is a Sylow $q$-subgroup of $G$.

Recall that $|\mathrm{IBr}_{p}^{1}(\overline{G})|=1$ and that $\mathrm{O}_{p}(\overline{G})=1$.
So, $\overline{G}$ is not nilpotent.
In fact, otherwise, $\overline{G}=G/\Phi$ is abelian, and hence $|\mathrm{IBr}_{p}^{1}(\overline{G})|=0$, a contradiction.
Therefore, part (3) of Lemma \ref{lem: ibrG/N} yields the first statement of Claim 2, while
\cite[Theorem A]{dolfi} yields the second.
Now, we assert that $F$ is $G$-indecomposable.
Otherwise, as $\overline{F}$ and $\Phi$ are two $G$-chief factors in $F$, we set $F=\Phi\times R$ where $R\unlhd G$.
Recall that $G$ acts both transitively on $\Phi^\sharp$ and $\overline{F}^\sharp$,
and hence there exist $\alpha \in \mathrm{Irr}(\Phi)^\sharp$, and $\beta \in \mathrm{Irr}(R)^\sharp$
such that 
$\mathrm{I}_{G}(\beta)<G$.
Setting $\gamma=\alpha\times \beta$, we have that $\mathrm{I}_{G}(\gamma)=\mathrm{I}_{G}(\alpha)\cap \mathrm{I}_{G}(\beta)<G$.
However, irreducible $p$-Brauer characters of $G$ lying over $\beta$ or $\gamma$ are all nonlinear by Clifford's theorem \cite[Corollary 8.7]{navarrobook}, and $\varphi \in \mathrm{IBr}_{p}^{1}(G|\alpha)$.
Therefore, $|\mathrm{IBr}_{p}^{1}(G)|\geq 3$, a contradiction.
Thus, $F=\mathrm{F}(G)$ is a $G$-indecomposable nilpotent subgroup of $G$.
In particular, $F$ is a $q$-group where $q$ is a prime distinct from $p$.
Moreover, by the second statement of Claim 2, it is straightforward to verify that $(|\overline{G}/\overline{F}|,|\overline{F}|)=1$.
Consequently, $F\in \mathrm{Syl}_{q}(G)$ and thus $G=F  \rtimes H$ by the Schur-Zassenhaus theorem.

\textbf{Claim 3.} $\mathrm{IBr}_{p}(F)^\sharp$ is a union of two $G$-orbits, say $\mathrm{IBr}_{p}(\overline{F})^\sharp$ and $\mathrm{IBr}_{p}(F)-\mathrm{IBr}_{p}(\overline{F})$;
$\mathrm{I}_{G}(\alpha)/F$ is a $p$-group for each $\alpha \in \mathrm{IBr}_{p}(F)^\sharp$.

As $\overline{H}$ acts transitively on $\overline{F}^\sharp$ by Claim 2,
it also acts transitively on $\mathrm{Irr}(\overline{F})^\sharp$ by \cite[Corollary 6.33]{isaacsbook}.
Let $\alpha \in \mathrm{IBr}_{p}(\overline{F})^\sharp$,
 and note that $\mathrm{IBr}_{p}(\overline{F})=\mathrm{Irr}(\overline{F})$. 
 Therefore,
by Claim 2 and Lemma \ref{lem: stabdual},
$\mathrm{I}_{\overline{G}}(\alpha)/\overline{F}$ is a $p$-group.
Since $\mathrm{I}_{G}(\alpha)/F\cong \mathrm{I}_{\overline{G}}(\alpha)/\overline{F}$,
$\mathrm{I}_{G}(\alpha)/F$ is also a $p$-group.
Take $\beta \in \mathrm{IBr}_{p}(F)-\mathrm{IBr}_{p}(\overline{F})= \mathrm{Irr}(F)-\mathrm{Irr}(\overline{F})$.
As $\Phi$ is minimal normal in $G$ by Claim 1, it is central in $F=\mathrm{F}(G)$, and thus,
$\beta_\Phi=\beta(1)\lambda$ for some $\lambda\in \mathrm{IBr}_p(\Phi)^\sharp$.
Recall that $\mathrm{IBr}_{p}(G|\lambda)=\{ \varphi \}$ and that $\varnothing\neq \mathrm{IBr}_{p}(G|\beta) \subseteq \mathrm{IBr}_{p}(G|\lambda)$ by \cite[Corollary 8.3]{navarrobook},
and hence $\mathrm{IBr}_{p}(G|\beta)=\{ \varphi \}$ for all $\beta \in \mathrm{IBr}_{p}(F)-\mathrm{IBr}_{p}(\overline{F})$.
Therefore, Clifford's theorem \cite[Corollary 8.7]{navarrobook} implies that $\mathrm{IBr}_{p}(F)-\mathrm{IBr}_{p}(\overline{F})$
forms a complete $G$-orbit.
Notably, $(|\mathrm{I}_{G}(\beta)/F|,|F|)=1$, implying that $\beta$ extends to $\mathrm{I}_{G}(\beta)$ by \cite[Theorem 8.13]{navarrobook}.
As $|\mathrm{IBr}_{p}(G|\beta)|=1$, an application of part (1) of Lemma \ref{lem: linearext} to $\mathrm{IBr}_{p}(G|\beta)$ yields that $|\mathrm{IBr}_{p}(\mathrm{I}_{G}(\beta)/F)|=1$.
Consequently,
$\mathrm{I}_{G}(\beta)/F$ is a $p$-group.

\smallskip

In the following, we set characters $\alpha \in \mathrm{IBr}_{p}(\overline{F})^\sharp$ and $\beta \in \mathrm{IBr}_{p}(F)-\mathrm{IBr}_{p}(\overline{F})$,
we denote $\mathrm{I}_{H}(\alpha)$ by $A$ and $\mathrm{I}_{H}(\beta)$ by $B$,
we write $|\Phi|=q^{d}$, $|\overline{F}|=q^{n}$, $a=|\mathrm{IBr}_{p}(\overline{F})^\sharp|$ and $b=|\mathrm{IBr}_{p}(F)-\mathrm{IBr}_{p}(\overline{F})|$.
By Claim 3, $A$ and $B$ are $p$-groups having indices $a$ and $b$ in $H$, respectively.

\smallskip

\textbf{Claim 4.} $F$ is a special $q$-group with center $\Phi$; 
$a=q^{n}-1$, $b=q^d-1$ and $a/b=p^k$ where $d$ is a proper divisor of $n$.

We begin by demonstrating that $F$ is nonabelian.
Working by contradiction, we assume that $F$ is abelian.
Hence, it follows that $|H:B|=|F|-|\overline{F}|=q^{n+d}-q^{n}$.
Note that $|H:A|=|\overline{F}|-1=q^{n}-1$, and hence
it yields that $q^{n}\mid |A|$, whereas $A$ is a $p$-group, a contradiction.
Therefore, $F$ is nonabelian.
Furthermore, as Claims 1 and 2 indicate, $\overline{F}$ and $\Phi$ are the $G$-chief factors in the nonabelian $q$-group $F$,
so $\Phi=\Phi(F)=F'=\mathrm{Z}(F)$.
In particular, $F$ is a special $q$-group.
As $\mathrm{IBr}_{p}(F)-\mathrm{IBr}_{p}(\overline{F})=\mathrm{Irr}(F)-\mathrm{Irr}(\overline{F})$,
it follows that $\beta(1)^2\mid |F:\mathrm{Z}(F)|=q^n$.
Write $\beta(1)=q^s$.
Counting the degree of the regular character of $F$, 
we have $q^n+b\cdot q^{2s}=q^{n+d}$.
By calculation, $b=q^{n-2s}(q^{d}-1)$.
As $|A|=|H|/(q^n-1)$ and $|B|=|H|/q^{n-2s}(q^{d}-1)$, we have the equation:
$$q^{n-2s}|B|=|A|\cdot \frac{q^{n}-1}{q^{d}-1}.$$
However, by Claim 3, both $A$ and $B$ are $p$-groups where $p\neq q$.
This implies that $n=2s$, $b=q^{d}-1$ and $a/b=(q^{n}-1)/(q^{d}-1)=p^{k}$.
So, we conclude that $d\mid n$.
Let $x\in F-\Phi$, $x^F$ an $F$-class containing $x$, and $t$ the number of $F$-classes in $F-\Phi$. 
According to Claim 3 and \cite[Theorem 13.24]{isaacsbook}, the set of nontrivial $G$-classes in $F$ consists of only two elements. 
As $G$ acts transitively on $\Phi^\sharp$ and $\Phi=\mathrm{Z}(F)$,
$\Phi^\sharp$ forms a complete $G$-class consisting of $q^d-1$ $F$-classes.
Consequently, $F-\Phi$ represents the other $G$-class in $F$, indicating that $|F-\Phi|=|x^F|t$.
The class equation of $F$ yields that
$q^d+|x^F|t=q^{n+d}$.
Note that $t=q^{n}+b-q^d=q^{n}-1$, and so
$|F:\mathrm{C}_{F}(x)|=|x^{F}|=q^{d}$.
In particular, $|\mathrm{C}_{F}(x)|=q^{n}$.
As $\mathrm{C}_{F}(x)>\mathrm{Z}(F)=\Phi$, we conclude that $n>d$.

\smallskip

\textbf{Claim 5.} Final conclusion.

Observe by Claim 4 that $F$ is a special $q$-group such that $H$ acts transitively on $\Phi^\sharp$ and that $q^{n}-1=p^{k}(q^{d}-1)$.
Assume that $\overline{G}\cong\mathsf{PrimitiveGroup}(q^n,i)$ where $(q^{n},i)\in \{ (5^{2},18),(3^{4},71) \}$, and $p=2$.
As $d$ is a proper divisor of $n$, a simple calculation implies that $|H:B|=q^d-1$ is a power of $2$.
However, $B$ is itself a $2$-group, which implies that $\overline{H} \cong H$ is also a 2-group which contradicts $q^{n}-1\mid |H|$.
Now, assume that $\overline{G}\cong A\Gamma(q^{2})$ where $q \in  \mathfrak{M}$ and $p=2$.
In this case, $H\cong \overline{H}\cong \Gamma(q^{2})$ acts transitively on both $\overline{F}^\sharp$ and $\Phi^\sharp$.
Since $\overline{F}\cong (\mathsf{C}_{q})^{2}$ and $\Phi\cong \mathsf{C}_{q}$, we deduce by Lemma \ref{lem: omega} that $F\cong \mathsf{ES}(q^{3}_+)$.
It is known that $\mathrm{Aut}(F)=\mathrm{Inn}(F) \rtimes \mathrm{GL}_2(q)$ by \cite[Theorem 1]{winter1972}.
Since $H$ acts faithfully on $F$, $H$ is isomorphic to a subgroup of $\mathrm{Aut}(F)$.
Note that 
$(|H|,|F|)=1$ as stated in Claim 2
and that the first cohomology group $H^{1}(\mathrm{GL}_2(q), \mathrm{Inn}(F))=0$ (see, for instance, \cite[Kapitel I, Aufgaben 62)]{huppertbook1}). 
Hence, up to a suitable conjugation in $\mathrm{Aut}(F)$,
we may identify $H$ as a subgroup of $\mathrm{GL}_2(q)$. 
Consequently, we conclude that $G\cong \mathsf{ES}(q^{3}_+) \rtimes \Gamma(q^{2})$ where $q \in  \mathfrak{M}$ and $p=2$. 
The remaining assertions in (1), (2) and (3) are supported by Claim 2 and Lemma \ref{lem: structure}.
 \end{proof}

\section{Nonsolvable groups}

In this section, we endeavor to classify finite nonsolvable groups with two nonlinear irreducible Brauer characters.
To achieve that, we need the classification of \emph{simple $K_3$-groups} (nonabelian simple groups of order divisible by exactly three distinct prime divisors).

Recall that $\mathrm{cd}(G) = \{ \chi(1):\chi \in \mathrm{Irr}(G) \}$ is the set of irreducible character degrees.

\begin{thm}\label{nonsol=2}
   Let $G$ be a nonsolvable group with $\mathrm{O}_{p}(G)=1$.
   Then $|\mathrm{IBr}_{p}^{1}(G)|=2$ if and only if $G$ is an almost simple group, and
   \[
     (G,p)\in \{ (\mathsf{A}_5,5),(\mathsf{S}_5,2), (\mathrm{PGL}_2 (7),2), (\mathrm{M}_{10},2), (\mathrm{Aut}(\mathsf{A}_6 ),2),  (\mathrm{P}\Sigma \mathrm{L}_2 (8),3)\}.
   \]
\end{thm}
\begin{proof}
   If $(G,p)\in \{ (\mathsf{A}_5,5),(\mathsf{S}_5,2), (\mathrm{PGL}_2 (7),2), (\mathrm{M}_{10},2), (\mathrm{Aut}(\mathsf{A}_6 ),2),  (\mathrm{P}\Sigma \mathrm{L}_2 (8),3)\}$, then $|\mathrm{IBr}_{p}^{1}(G)|=2$ by checking via $\mathsf{GAP}$ \cite{gap}.

   Conversely, assume that $|\mathrm{IBr}_{p}^{1}(G)|=2$.
   Then $p\mid |G|$.
   In fact, if $p\nmid |G|$, then $\mathrm{IBr}_{p}(G)=\mathrm{Irr}(G)$ and $|\mathrm{cd}(G)|\leq 3$, and hence $G$ is solvable by \cite[Corollary 12.6, Theorem 12.15]{isaacsbook}.
  We claim now that $\mathrm{F}(G)=1$. In fact, otherwise $|\mathrm{IBr}_{p}^{1}(G/\mathrm{F}(G))|\leq 1$ by part (1) of Lemma \ref{lem: ibrG/N}, and so $G/\mathrm{F}(G)$ is solvable by the observations proceeding Lemma \ref{lem: transitive action of G on V-0}, whereas $G$ is nonsolvable.
   Write $S=\mathrm{F}^*(G)$ where $\mathrm{F}^*(G)$ denotes the generalized Fitting subgroup of $G$.
   Since $\mathrm{F}(G)=1$,  
   the generalized Fitting subgroup $S$ is a direct product of some nonabelian simple groups.
   Notably, $S$ has a unique linear $p$-Brauer character.
   Next, we show that $\mathrm{IBr}_{p}^{1}(S)$ consists of at most two $G$-orbits.
  In fact, for $\theta\in\mathrm{IBr}_{p}(S)^\sharp$, $\varphi(1)>1$ for each $\varphi \in \mathrm{IBr}_{p}(G|\theta)$ by Clifford's theorem \cite[Corollary 8.7]{navarrobook};
   as $|\mathrm{IBr}_{p}^{1}(G)|=2$, $\mathrm{IBr}_{p}(S)^\sharp$ is a union of at most two $G$-orbits.
   Consequently, $\mathrm{IBr}_{p}(S)$ contains at most three $G$-orbits, implying that 
   $|\omega_{p'}(S)|\leq 3$. 
  Recall that $S=\mathrm{F}^*(G)$ is a direct product of some nonabelian simple groups,
  and hence $S$ is a simple group.
  In fact, otherwise, $S=S_1\times \cdots \times S_t$ where $S_i$ are nonabelian simple groups and $t>1$;
  as $|\pi(S_i)|\geq 3$ by Burnside's $p^{a}q^{b}$-theorem, there exist $x_i,y_i\in S_i$ such that $o(x_i)=r_i$ and $o(y_i)=\ell_i$ where $r_i,\ell_i\in \pi(S_i)-\{ p \}$;
  so, $\{ 1, r_1, \ell_1, r_1r_2 \} \subseteq \omega_{p'}(S)$ which contradicts $|\omega_{p'}(S)|\leq 3$. 
 Further, as $|\pi(S)|\leq |\omega_{p'}(S)|\leq 3$, $S$ is a simple $K_3$-group.
  Since the classification of simple $K_3$-groups
  yields that 
  $$S\in \{ \mathsf{A}_5,  \mathrm{PSL}_2(7), \mathsf{A}_6, \mathrm{PSL}_2(8), \mathrm{PSL}_2(17), \mathrm{PSL}_3(3), \mathrm{PSU}_3(3), \mathrm{PSU}_4(2) \}$$ (see, for instance, \cite[Page 12]{gorensteinbook}),
  our desired results follow from direct calculations using $\mathsf{GAP}$ \cite{gap}.
\end{proof}

\section{Proof of Theorem \ref{thmA}}

We begin with three preliminary lemmas.

\begin{lem}\label{wreathp}
  Let $p$ be an odd prime, and $G=A \wr Y$ where $A=\langle a\rangle$ has order $2p^n$
  and $Y=\langle y\rangle$ has order $2$.
  If $H$ is a maximal subgroup of $G$ isomorphic to $\mathsf{Dic}_{p^n}\times \mathsf{C}_{p^n}$,
  then $H=\langle ay_0, a^{2}\rangle$ for some involution $y_0\in G-AA^y$.
\end{lem}
\begin{proof}
  Let $P \in \mathrm{Syl}_{p}(G)$, $Q \in \mathrm{Syl}_{2}(G)$ and $B=AA^y$.
  Then $P=\langle a^2,a^{2y}\rangle \unlhd G$.
  Also, the Schur-Zassenhaus theorem yields that $G=P \rtimes Q$ where $Q\cong G/P\cong  \mathsf{D}_8$.
  Note that $H\cong\mathsf{Dic}_{p^n}\times \mathsf{C}_{p^n}$, and so $H=P  \rtimes (Q\cap H)$ where $Q\cap H\cong \mathsf{C}_{4}$.
  We next show that $Q\cap H=\langle a^{p^n}y_0\rangle$ for some involution $y_0\in G-B$.
  Since $Q\cap B=\langle a^{p^{n}},a^{p^{n}y}\rangle\cong (\mathsf{C}_{2})^{2}$ and $Q\cong \mathsf{D}_8$,
  we take an involution $y_0\in Q-(Q\cap B)\subseteq G-B$.
  So, $\langle a^{p^n}y_0\rangle$ is the unique cyclic subgroup of order $4$ in $Q$.
  Consequently, $Q\cap H=\langle a^{p^n}y_0\rangle$.
  Observe that $2k+p^n=1$ where $k=(1-p^n)/2\in \mathbb{Z}$ and that $a^2, a^{p^n}y_0 \in H$,
  and so $ay_0=(a^{2})^ka^{p^n}y_0\in H$.
 Thus, we establish that $H=\langle ay_0,a^{2}\rangle$.
\end{proof}

\begin{lem}\label{imprim}
  Let $H$ be a finite group and $V=V_1\times V_2$ an imprimitive $H$-module over a prime field $\mathbb{F}_q$.
  Assume that $H\leq T\wr \langle \sigma\rangle$ where $T=\Gamma_0(V_1)=\langle t\rangle$ and $\sigma$ is an involution transposing $V_1$ and $V_2$.
  If $t\sigma, t^2\in H$, then $V^\sharp$
  is a union of two distinct $H$-orbits $V_1^\sharp \sqcup V_2^\sharp$ and $V-(V_1\cup V_2)$.
  In addition, if either $H<T\wr \langle \sigma\rangle$ and $|V_1|=2p^s+1$ for an odd prime $p$, or $|V_1|=p^{s}+1$ where $p=2$, then $H'$ and $\mathrm{C}_{H}(v)$ are $p$-groups for all $v \in V^\sharp$.
\end{lem}
\begin{proof}   
  Since $t\sigma, t^{2}\in H$, it follows that $t t^{\sigma}=(t\sigma)^{2}$ and $t^{2 \sigma}=t^{-2}(t\sigma)^{4}$ both lie in $H$.
  Let $v_1\in V_1^\sharp$, and set $v=v_1$, $w=v_1v_1^{\sigma}$.
  Then $v$ and $w$ lie in distinct $H$-orbits.
  Observe that $v^{(t \sigma)^{2k}}=v_1^{t^{k}}$ 
  and $v^{(t \sigma)^{2k+1}}=(v_1^{t^{k+1}})^{\sigma}$ for $k\in \mathbb{Z}$,
  and that $T$ acts transitively on $V_1^\sharp$.
  As $v$ and $u_1u_2$ are not conjugate in $H$ for each $u_i\in V_i^\sharp$,
  it follows that $v^H=V_1^\sharp \sqcup V_2^\sharp$.
  Also, note that $w^{t^{2k}t\sigma}=v_1(v_1^{t^{2k+1}})^{\sigma}$ 
  and $w^{t^{2k }\sigma}=v_1 (v_1^{t^{2k}})^{\sigma}$,
  and hence $\{ v_1 u_1^\sigma: u_1\in V_1^\sharp \}\subseteq w^{H}$.
  As $(v_1 u_1^\sigma)^{t \sigma}=u_1  (v_1^t)^\sigma$ and
  $(u_1^{t^{-k}} (v_1^t)^\sigma)^{(t\sigma)^{2k}}=u_1 (v_1^{t^{k+1}})^\sigma$ for $k\in \mathbb{Z}$,
  we conclude that $w^{H}=V-(V_1\cup V_2)$.

  If in addition $|V_1|=p^s+1$ where $p=2$, then $T\wr \langle \sigma\rangle$ is a $2$-group, and so our desired result follows.
  Assume in addition that
   $H<T\wr \langle \sigma\rangle$ and $|V_1|=2p^s+1$ for an odd prime $p$.
   Write $X=\langle t^2, tt^\sigma\rangle \langle t\sigma\rangle$.
   As $\langle t^2\rangle\cap \langle tt^\sigma\rangle=1$ and $t\sigma\not \in \langle t^2, tt^\sigma\rangle$,
   $|T\wr \langle \sigma\rangle:X|\leq 2$.
   Since $X\leq H<T\wr \langle \sigma\rangle$, $H=X$ has index 2 in $T\wr \langle \sigma\rangle$.
   In particular, $|H|=(|V_1|-1)^{2}=4p^{2s}$.
   So, $\mathrm{C}_{H}(v)$ is a $p$-group for each $v\in V^\sharp$ by calculation.
   As $H=\langle t^2, tt^\sigma\rangle \langle t\sigma\rangle$ where $\langle t^2, tt^\sigma\rangle=\langle t^2\rangle\times \langle tt^\sigma\rangle$, we conclude that $H'=[\langle t^2, tt^\sigma\rangle, \langle t\sigma\rangle]=\langle [t^{2},\sigma]\rangle\cong \mathsf{C}_{p^s}$.
\end{proof}

\begin{lem}\label{ibrG=2}
   Let $G$ be a solvable primitive group of rank $3$ with $\mathrm{O}_{p}(G)=1$.
   Write $G=V \rtimes H$ where $V=\mathrm{F}(G)$ and $H$ a complement of $V$ in $G$.
   If $H'$ and $\mathrm{C}_{H}(v)$, where $v \in V^\sharp$, are $p$-groups, then $|\mathrm{IBr}_{p}^{1}(G)|=2$.
\end{lem}
\begin{proof}
   Since $H'$ is a $p$-group, $H'\leq \mathrm{O}_{p}(H)$, and so $|\mathrm{IBr}_{p}^{1}(G/V)|=|\mathrm{IBr}_{p}^{1}(H)|=|\mathrm{IBr}_{p}^{1}(H/\mathrm{O}_{p}(H))|=0$.
   As $V$ is a $p'$-group, $\mathrm{IBr}_{p}(V)=\mathrm{Irr}(V)$.
   Furthermore, since $V^\sharp$ is a union of precisely two $G$-orbits, we deduce that $\mathrm{Irr}(V)^\sharp$ is also a union of exactly two $G$-orbits by \cite[Corollary 6.33]{isaacsbook}, say $\mathcal{O}_\lambda$ and $\mathcal{O}_\mu$ where $\lambda,\mu \in \mathrm{Irr}(V)^\sharp$.
   Note that $V$ is $H$-irreducible if and only if $\mathrm{Irr}(V)$ is $H$-irreducible,
   and so both $\mathcal{O}_\lambda$ and $\mathcal{O}_\mu$ have size larger than 1.
   Therefore, by Clifford's theorem \cite[Corollary 8.7]{navarrobook}, we conclude that $\mathrm{IBr}_{p}(G|\alpha)=\mathrm{IBr}_{p}^1(G|\alpha)$ for every $\alpha \in \mathrm{Irr}(V)^\sharp$.
   For each $\alpha \in \mathrm{Irr}(V)^\sharp$, noting that $G=V \rtimes H$ and by Lemma \ref{lem: stabdual} that $\mathrm{I}_{H}(\alpha)$ is a $p$-group, 
   we establish that $|\mathrm{IBr}_{p}(G|\alpha)|=|\mathrm{IBr}_{p}(\mathrm{I}_{H}(\alpha))|=1$ by part (2) Lemma \ref{lem: linearext}.
   As $\mathrm{IBr}_{p}^{1}(G)=\mathrm{IBr}_{p}^{1}(G/V)\sqcup \mathrm{IBr}_{p}(G|\lambda)\sqcup \mathrm{IBr}_{p}(G|\mu)$, this leads to the conclusion that $|\mathrm{IBr}_{p}^{1}(G)|=2$.
\end{proof}

Now, we are ready to prove Theorem \ref{thmA}.
In the proof, we use $\mathrm{Irr}_1(G)$ to denote the set of nonlinear irreducible characters of $G$.

\begin{proof}[Proof of Theorem \ref{thmA}]
  By Theorem \ref{nonsol=2}, we may assume that $G$ is solvable.
 Let $F=\mathrm{F}(G)$, $\Phi=\Phi(G)$ and $\overline{G}=G/\Phi$.
  By Gasch\"{u}tz's theorem, $\overline{G}=\overline{F}  \rtimes \overline{H}$ where $\overline{F}$ is a completely reducible and faithful $\overline{H}$-module (possibly of mixed
 characteristic).
 As $\mathrm{O}_{p}(G)=1$, $F=\mathrm{F}(G)$ is a $p'$-group.
 Hence, $\mathrm{F}(\overline{G})=\overline{F}$ is also a $p'$-group.

 \smallskip

  Assume that $|\mathrm{IBr}_{p}^{1}(G)|=2$.
  Then $G$ is nonabelian.
  If $|\mathrm{IBr}_{p}^{1}(\overline{G})|=1$, then either (5) or (6) holds by Proposition \ref{fnonabel}.

  \textbf{Claim A}. If $|\mathrm{IBr}_{p}^{1}(\overline{G})|=0$, then $(1)$ holds.

  As $|\mathrm{IBr}_{p}^{1}(\overline{G})|=0$, it follows that $\overline{G}=\overline{G}/\mathrm{O}_{p}(\overline{G})$ is abelian, indicating that $G$ is nilpotent.
  Given that $\mathrm{O}_{p}(G)=1$, $G$ is a $p'$-group.
  In particular, $\mathrm{IBr}_{p}(G)=\mathrm{Irr}(G)$.

  Let $1<Z\unlhd\mathrm{Z}(G)$.
  We claim first that if $G/Z$ is nonabelian then $G/Z$ is an extraspecial $2$-group and $Z\cong \mathsf{C}_{2}$.
  Since $Z>1$,
  it follows by part (1) of Lemma \ref{lem: ibrG/N} that $|\mathrm{Irr}_1(G)-\mathrm{Irr}_1(G/Z)|=|\mathrm{Irr}_1(G/Z)|=1$.
 As $G/Z$ is a $p'$-group,
  our conclusion follows directly by \cite[Theorem A]{dolfi} and part (2) of Lemma \ref{lem: ibrG/N}.

  We claim next that $G'$ is minimal normal in $G$.
  Assume not and let $Z$ be a minimal $G$-invariant subgroup of $G'$.
  As $G/Z$ is nonabelian, $G/Z$ is an extraspecial $2$-group, $Z\cong \mathsf{C}_{2}$ and $|G'|=4$.
  Write $\mathrm{Irr}_{1}(G/Z)=\{ \varphi \}$, $\mathrm{Irr}_{1}(G)-\mathrm{Irr}_{1}(G/Z)=\{ \chi \}$
  and $|G/G'|=2^{2n}$.
  Then $\varphi(1)^2=2^{2n}$.
  However, since $\chi(1)^{2}+\varphi(1)^2+|G/G'|=|G|$, we have $\chi(1)^{2}=2^{2n+1}$, a contradiction.

  Set $\mathrm{Irr}_1(G)=\{ \chi, \varphi \}$ and $|G'|=q$ where $q$ is a prime.
  As $G'$ is minimal normal in $G$, it follows by \cite[Problem 2.13]{isaacsbook} that $\chi(1)^{2}=\varphi(1)^{2}=|G/\mathrm{Z}(G)|$.
Note that $\chi(1)^{2}+\varphi(1)^{2}+|G/G'|=|G|$ and that $G'\leq \mathrm{Z}(G)$,
  and hence $|\mathrm{Z}(G)|=2q/(q-1)$.
 By calculation, $q=2$ or $3$.
 If $q=3$, then $G'=\mathrm{Z}(G)$ is the unique minimal normal subgroup of $G$, and hence $G$ is an extraspecial $3$-group by \cite[Lemma 12.3]{isaacsbook}.
 Assume now that $q=2$.
 This implies $|\mathrm{Z}(G)|=4$, and so $G$ is a $2$-group.

  \textbf{Claim B.} If $|\mathrm{IBr}_{p}^{1}(\overline{G})|=2$, then 
  one of (2), (3) or (4) holds.

 As $|\mathrm{IBr}_{p}^{1}(\overline{G})|=|\mathrm{IBr}_{p}^{1}(G)|$ and $G$ is nonabelian,
 part (1) of Lemma \ref{lem: ibrG/N} implies that $\Phi=1$.
 Consequently, $G=F \rtimes H$ where $F=\mathrm{F}(G)$ is a completely reducible and faithful $H$-module (possibly of mixed
 characteristic).

 Let $Z=\mathrm{Z}(G)$.
 As $Z$ is an $H$-submodule of $F$, there exists an $H$-submodule $V$ of $F$ such that $F=V\times Z$.
 Since $G$ is nonnilpotent, we have $V>1$.
 Given that $V$ is abelian and $V\cap Z=1$, we claim that $\mathrm{I}_{H}(\mu)<H$ for all $\mu \in \mathrm{Irr}(V)^\sharp$.
 In fact, take an $H$-invariant $\mu \in \mathrm{Irr}(V)$, $\ker(\mu)$ is an $H$-submodule of $V$ and $[V,H]\leq \ker(\mu)$; as $V$ is a completely reducible $H$-module, $V=\ker(\mu)\times W$ where $W$ is an $H$-submodule of $V$;
 hence, $[W,H]\leq W\cap \ker(\mu)=1$, i.e. $W$ is central in $G$;
 so, $W=1$ and $V=\ker(\mu)$, i.e. $\mu=1_V$.
 For $\mu\in\mathrm{Irr}(V)^\sharp$ and $\zeta\in \mathrm{Irr}(Z)$, as $\mathrm{I}_{H}(\mu\times \zeta)=\mathrm{I}_{H}(\mu)\cap \mathrm{I}_{H}(\zeta)<H$,
 Clifford's theorem \cite[Corollary 8.7]{navarrobook} yields that $\varphi(1)>1$ for all $\varphi \in \mathrm{IBr}_{p}(G|\mu\times \zeta)$.
 Note that $|\mathrm{IBr}_{p}^{1}(G)|=2$, and hence either $Z=1$ and $\mathrm{Irr}(V)^\sharp$ contains at most two $H$-orbits or $|Z|=2$.

 Assume that $|Z|=2$.
 Then $p>2$ and $G=(VH)\times Z$ where $Z\cong \mathsf{C}_{2}$.
 Since $\mathrm{IBr}_{p}(G)=\{ \alpha \times \zeta: \alpha \in \mathrm{IBr}_{p}(VH),\zeta\in \mathrm{IBr}_{p}(Z) \}$ by \cite[Theorem 8.21]{navarrobook} and $|\mathrm{IBr}_{p}^{1}(G)|=2$,
 it follows that $|\mathrm{IBr}_{p}^{1}(VH)|=1$.
 Note that $\mathrm{O}_{p}(VH)\leq \mathrm{O}_{p}(G)=1$ and that $p>2$,
 and hence \cite[Theorem A]{dolfi} yields that $VH$ is a doubly transitive Frobenius group with complement $H$ such that $H'$ is a $p$-group, i.e. (4) holds.

 If $Z=1$ and $\mathrm{Irr}(V)^\sharp$ contains at most two $H$-orbits,
 then $V=F$.
 According to \cite[Corollary 6.33]{isaacsbook},
 the number of $H$-orbits in $\mathrm{Irr}(V)^\sharp$ equals the number of $H$-orbits in $V^\sharp$.
 Since $\mathrm{Irr}(V)^\sharp$ contains at most two $H$-orbits,
 $V^\sharp$ also contains at most two $H$-orbits.
  So, $V$ is $H$-irreducible.
 Write $|V|=q^n$ where $q$ is a prime distinct from $p$.
 As $V=F=\mathrm{F}(G)$, $V$ is the unique minimal normal subgroup of $G$.
 Consequently, $G$ is a solvable primitive permutation group.
 In particular, if $V^\sharp$ is a full $H$-orbit, then $G$ is a solvable doubly transitive permutation group; if $V^\sharp$ contains two $H$-orbits, then $G$ is a
 solvable primitive permutation group of rank 3.
 Therefore, Propositions \ref{prop: rank 2 1}, \ref{prop: rank 2 2} and \ref{prop: rank3 1} imply (2),
 and Proposition \ref{prop: rank3 2} implies (3).
 
\smallskip

 Conversely, suppose that $G=F  \rtimes H$, where $F=\mathrm{F}(G)$, is a group listed in one of (1), (2), (3), (4), (5) or (6).
 Then $F$ is a $p'$-group. 
 In particular, $\mathrm{IBr}_{p}(F)=\mathrm{Irr}(F)$.
 If $p=2$, and $G\cong\mathsf{PrimitiveGroup}(q^{n},i)$ where
 $
  (q^{n},i)\in \{ (5^{2},19),(7^{2},25),(3^{4},99),(7^{2},18),(17^{2},82),(23^{2},51) \}
 $,
 then $|\mathrm{IBr}_{p}^{1}(G)|=2$ by calculation via $\mathsf{GAP}$ \cite{gap}.
 If one of (2a), (2b), (2f) or (3b) holds, 
 then $|\mathrm{IBr}_{p}^{1}(G)|=2$ by calculation via $\mathsf{GAP}$ \cite{gap}.
 If (2d) holds, then $G$ has exactly two nonlinear irreducible $p$-Brauer characters of the same degree $(q^n-1)/2$.

 If (1) holds,
 then $p\nmid |G|$ and $G$ is a nonabelian $q$-group where $q\in \{ 2,3 \}$ such that $G'$ is a minimal normal subgroup of $G$.
 So, $\mathrm{IBr}_{p}(G)=\mathrm{Irr}(G)$ and $G'\leq \mathrm{Z}(G)$.
 Also, by \cite[Problem 2.13]{isaacsbook}, $\chi(1)^{2}=|G:\mathrm{Z}(G)|$ for all nonlinear irreducible characters $\chi\in \mathrm{Irr}(G)$.
 Let $t=|\mathrm{Irr}_1(G)|$.
 Write $|G:\mathrm{Z}(G)|=q^{2n}$ and $|\mathrm{Z}(G)|=q^k$.
 Then $k=2$ if $q=2$ and $k=1$ if $q=3$.
 As $|G|=\sum_{\chi \in \mathrm{Irr}(G)}\chi(1)^2$, $q^{2n+k}=q^{2n+k-1}+t q^{2n}$.
 It is straightforward to check that $t=2$.

 Assume that (2c) holds.
 Then $p>2$, $G=M \rtimes C$, $H=L \rtimes C$ and, $M=F \rtimes L$ is a doubly transitive Frobenius group with kernel $F$ and complement $L$.
 Also, $H'$ is a $p$-group and $C=\mathrm{C}_{H}(x_0)\cong \mathsf{C}_{2}$ for some $x_0\in F^\sharp$.
 As $H'$ is a $p$-group, $L'$ is also a $p$-group.
 Note that $M$ is a doubly transitive Frobenius group with kernel $F$ and complement $L$ such that $L'$ is a $p$-group 
 and that $\mathrm{O}_{p}(M)=1$.
 So, \cite[Theorem A]{dolfi} yields that $\mathrm{IBr}_{p}^{1}(M)=\{ \theta \}$ where $\theta=\alpha^M$ for each $\alpha \in \mathrm{Irr}(F)^\sharp$.
 Observing that the cyclic group $C=\mathrm{C}_{H}(x_0)$ fixes some $\lambda \in \mathrm{Irr}(F)^\sharp$,
 we deduce that $\theta=\lambda^{M}$ is $G$-invariant. 
As $G/M\cong \mathsf{C}_{2}$,
$\theta$ extends to $\hat{\theta}\in \mathrm{IBr}_{p}(G)$ by \cite[Theorem 8.12]{navarrobook}.
 So, by \cite[Corollary 8.20]{navarrobook}, $\mathrm{IBr}_{p}(G|\lambda)=\{ \hat{\theta}, \hat{\theta}\beta \}$ where $\beta\in \mathrm{IBr}_{p}(C)^\sharp$.
 As $H'$ is a $p$-group, $|\mathrm{IBr}_{p}^{1}(G/F)|=|\mathrm{IBr}_{p}^{1}(H)|=0$.
 Since $H$ acts transitively on $\mathrm{Irr}(F)^\sharp$, we conclude that $|\mathrm{IBr}_{p}^{1}(G)|=2$.

 \textbf{Claim 1.} If one of (2e), (2g), (3a) or (3c) holds, then $|\mathrm{IBr}_{p}^{1}(G)|=2$.

 By Lemma \ref{ibrG=2}, it suffices to show that $G$ is a solvable primitive permutation group of rank 3 such that $H'$ and $\mathrm{C}_{H}(x)$, where $x \in F^\sharp$, are $p$-groups.
 Note that $F=\mathrm{F}(G)$ is a minimal normal subgroup of the solvable group $G$,
 and so $G$ is a solvable primitive permutation group.
 Therefore, it suffices to show that $F^\sharp$ consists of exactly two $G$-classes 
 and that $H'$ and $\mathrm{C}_{H}(x)$, where $x\in F^\sharp$, are $p$-groups.

Assume that either (3a) or (3c) holds. 
Then we are done by Lemmas \ref{lem: wreath2}, \ref{wreathp} and \ref{imprim}.

 Assume that (2e) holds.
 Then $p=2$, $q\in \mathfrak{M}$, $F\cong (\mathsf{C}_{q})^{2}$ and $H$ is a subgroup of $\Gamma(q^{2})$ which is isomorphic to $\mathsf{D}_{2(q+1)}\times \mathsf{C}_{(q-1)/2}$.
 So, $H'\cong \mathsf{C}_{(q+1)/2}$ is a $2$-group.
 Let $H_0$ be the maximal cyclic subgroup of $H$, and observe that $H_0\cong \mathsf{C}_{(q^{2}-1)/2}$.
 Therefore, $F$ is a faithful irreducible $H_0$-module, and so $H_0$ acts fixed-point-freely on $F$.
 Noting that $|H:H_0|=2$, we conclude that $|\mathrm{C}_{H}(x)|\leq 2$ and $\frac{q^{2}-1}{2}\mid |x^G|$ for each $x\in F^\sharp$.
 Also, since $H\cong \mathsf{D}_{2(q+1)}\times \mathsf{C}_{(q-1)/2}$, $H$ cannot act fixed-point-freely on $F$.
 Consequently, $F^\sharp$ consists of exactly two $G$-classes.

 Assume that (2g) holds.
 Then $p=2$, $q= 3\cdot 2^k-1\geq 11$, $F\cong (\mathsf{C}_{q})^{2}$ and $H$ is a subgroup of $\Gamma(q^2)$ which is isomorphic to $\mathsf{SD}_{4(q+1)/3}\times \mathsf{C}_{(q-1)/2}$.
 Also, there exists some $x_0\in F^\sharp$ such that $\mathrm{C}_{H}(x_0)=1$.
 So, $H'\cong \mathsf{C}_{(q+1)/3}$ is a $2$-group.
 Let $H_0$ be the maximal cyclic subgroup of $H$, and observe that $H_0\cong \mathsf{C}_{(q^{2}-1)/3}$.
 Therefore, $F$ is a faithful irreducible $H_0$-module, and so $H_0$ acts fixed-point-freely on $F$.
 Noting that $|H:H_0|=2$, we conclude that $|\mathrm{C}_{H}(x)|\leq 2$ and $\frac{q^{2}-1}{3}\mid |x^G|$ for each $x\in F^\sharp$.
 Since there exists some $x_0\in F^\sharp$ such that $\mathrm{C}_{H}(x_0)=1$,
 we conclude that $F^\sharp$ consists of exactly two $G$-classes.
 Indeed, as $H\cong \mathsf{SD}_{4(q+1)/3}\times \mathsf{C}_{(q-1)/2}$, $H$ cannot act fixed-point-freely on $F$.

 \textbf{Claim 2.} If one of (4), (5) or (6) holds, then $|\mathrm{IBr}_{p}^{1}(G)|=2$.

Assume that (4) holds.
Note that $VH$ is a doubly transitive Frobenius group with complement $H$ such that $H'$ is a $p$-group, and that $\mathrm{O}_{p}(VH)\leq \mathrm{O}_{p}(G)=1$.
So, \cite[Theorem A]{dolfi} yields that $|\mathrm{IBr}_{p}^{1}(VH)|=1$.
Observe that $p>2$ and $G=(VH)\times Z$ where $Z=\mathrm{Z}(G)\cong \mathsf{C}_{2}$,
and hence it follows by \cite[Theorem 8.21]{navarrobook} that $|\mathrm{IBr}_{p}^{1}(G)|=2$.

Assume that either (5) or (6) holds.
It is straightforward to check that $(|G/F|,|F|)=1$; $H/\mathrm{O}_{p}(H)$ is abelian; $H$ acts transitively on both $\overline{F}^\sharp$ and $\Phi^\sharp$;
and $\mathrm{C}_{H}(\overline{x})$ and $\mathrm{C}_{H}(y)$ are $p$-groups for all $\overline{x} \in \overline{F}^\sharp$ and $y \in \Phi^\sharp$.
Moreover, $F$ is a nonabelian $q$-group (where $q\neq p$), and
for every $\theta \in \mathrm{IBr}_{p}^{1}(F)$,
there exists a unique $\lambda \in \mathrm{IBr}_p(\Phi)^\sharp$ 
such that $\theta=\frac{1}{\theta(1)}\lambda^{F}$.
Since $H/\mathrm{O}_{p}(H)$ is abelian, $|\mathrm{IBr}_{p}^{1}(G/F)|=|\mathrm{IBr}_{p}^{1}(H)|=0$.
As $\mathrm{C}_{H}(\overline{x})$ and $\mathrm{C}_{H}(y)$ are $p$-groups for all $\overline{x} \in \overline{F}^\sharp$ and $y \in \Phi^\sharp$,
by Lemma \ref{lem: stabdual} $\mathrm{I}_{H}(\alpha)$ and $\mathrm{I}_{H}(\lambda)$ are $p$-groups for all $\alpha \in \mathrm{IBr}_{p}(\overline{F})^\sharp$ and $\lambda \in \mathrm{IBr}_{p}(\Phi)^\sharp$
(as $F$ is a $p'$-group, $\mathrm{IBr}_{p}(F)=\mathrm{Irr}(F)$).
Now, let $\beta\in \mathrm{IBr}_{p}^{1}(F)$.
Then $\beta=\frac{1}{\beta(1)}\mu^{F}$ for a unique $\mu\in \mathrm{IBr}_{p}(\Phi)^\sharp$,
so $\mathrm{I}_{H}(\beta)~(\leq \mathrm{I}_{H}(\mu))$ is a $p$-group.
As $(|G/F|,|F|)=1$, $\alpha$ extends to $\hat{\alpha}\in \mathrm{IBr}_{p}(\mathrm{I}_{G}(\alpha))$ for all $\alpha \in \mathrm{IBr}_{p}(F)^\sharp$ by \cite[Theorem 8.13]{navarrobook}.
By part (1) of Lemma \ref{lem: linearext}, $|\mathrm{IBr}_{p}(G|\alpha)|=|\mathrm{IBr}_{p}(\mathrm{I}_{G}(\alpha)/F)|=1$ for all $\alpha \in \mathrm{IBr}_{p}(F)^\sharp$. 
Assume that $G$ is not isomorphic to $\mathrm{SL}_2(3)$.
Then $\varphi(1)\geq |G:\mathrm{I}_{G}(\alpha)|>1$ for each $\varphi \in \mathrm{IBr}_{p}(G|\alpha)$ and $\alpha \in \mathrm{IBr}_{p}(F)^\sharp$.
Recall that $H$ acts transitively on $\overline{F}^\sharp$ and $\Phi^\sharp$ and that $(|H|,|F|)=1$,
and hence it follows that $\mathrm{IBr}_{p}(F)^\sharp$ is a union of exactly two $H$-orbits.
Consequently, $|\mathrm{IBr}_{p}^{1}(G)|=2$.
If $p=3$ and $G\cong \mathrm{SL}_2(3)$,
then $G$ has only two nonlinear $p$-Brauer characters such that one has degree 2 and the other has degree 3.
\end{proof}

\begin{acknowledgement}
	The authors are grateful to  the referees for her/his
 valuable comments.
\end{acknowledgement}

\end{document}